\documentclass[12pt]{article}
\usepackage{amssymb}
\usepackage{amsmath}
\usepackage{eucal}

\newtheorem{theo}{Theorem}[section]
\newtheorem{lemma}[theo]{Lemma}

\def \P  {\mathbb{P}}
\def \Q  {\mathbb{Q}} 
\def \R  {\mathbb{R}} 
\def \Z  {\mathbb{Z}}  
\def \N  {\mathbb{N}} 
\def \ind {\hbox{ 1\hskip -3pt I}}
\def \E {\mathbb{E}}

\def \AA {{\cal A}}
\def \BB {{\cal B}}
\def \CC {{\cal C}}
\def \DD {{\cal D}}
\def \II {{\cal I}}
\def \JJ {{\cal J}}
\def \LL {{\cal L}}
\def \MM {{\cal M}}
\def \NN {{\cal N}}
\def \TT {{\cal T}}
\def \VV {{\cal V}}

\newcommand {\refeq}[1] {(\ref{#1})}
\newcommand {\acc}[1] {\left\{ {#1} \right\}}
\newcommand {\bra}[1] {\left\langle {#1} \right\rangle}
\newcommand {\cro}[1] {\left[ {#1} \right]}
\newcommand {\pare}[1] {\left( {#1} \right)}
\newcommand {\nor}[1] { \left\| {#1} \right\|}
\newcommand {\va}[1] {\left| {#1} \right|}

\def \xir {\bar{\xi}_r}
\def \xin {\bar{\xi}_{r_n}}
\def \pd {\psi_{\delta}}

\begin{document}

\vspace*{.2cm}
\begin{center}
{\Large{\bf LARGE DEVIATIONS FOR BROWNIAN MOTION IN A RANDOM SCENERY.}}
\vspace{5mm}\\
 {\large{\bf Amine Asselah}} \footnote{E-mail: asselah@cmi.univ-mrs.fr},
 {\large{\bf Fabienne Castell}} \footnote{E-mail: castell@cmi.univ-mrs.fr}.\\
 \vspace{.5mm}
 Laboratoire d'Analyse, Topologie et Probabilit\'es. CNRS UMR 6632. \\
 CMI. Universit\'e de Provence. \\ 39 rue Joliot Curie.  \\
 13453 Marseille Cedex 13. FRANCE.
\end{center}
\vspace{1cm}

{\bf Abstract.}
{\small We prove large deviations principles in large time, for 
the Brownian occupation time in random scenery 
$\frac{1}{t} \int_0^t \xi(B_s) \, ds $. 
The random field is constant on the elements of a partition 
of $\R^d$ into unit cubes. These random constants, say 
$\acc{\xi(j), j \in \Z^d}$ consist of i.i.d. bounded 
variables, independent of the Brownian motion $\{B_s,s\ge 0\}$. 
This model is a time-continuous version of Kesten and Spitzer's random 
walk in random scenery. We prove large deviations principles in 
``quenched'' and ``annealed'' settings. }
\vspace{.5cm}


Large deviations. Additive functionals of Brownian motion.
Random media.

{\em Running title: } Brownian motion in random scenery.

\vspace{1cm}
\section{Introduction.}

We study the large time asymptotics of random additive
functionals of Brownian motion $\int_0^t \xi(B_s)ds$, where the
random field $\{\xi(x), x\in \R^d\}$ is independent of the Brownian
motion $\{B_s,s\ge 0\}$. We consider the case where $\xi$ is
a random constant, say $\xi(i)$, on the $i^{\tt th}$ cube
of a partition of $\R^d$ into unit cubes. The sequence
$\{\xi(i),\ i\in \Z^d\}$ consists of i.i.d. bounded random
variables with common law $\nu_i=\nu$, and 
we assume for convenience that $|\xi(0)|\le 1$ and $E_{\nu}[\xi(0)]=0$.

This is related to one of Kesten-Spitzer's models of random walk 
in random scenery: let $\{X_i,\ i\in \N\}$ be a sequence
of $\Z^d$-valued i.i.d. random vectors with mean 0 and finite non-singular
covariance matrix $\Sigma$, and define $S_n=X_1+\dots+X_n$.
Let $\{\xi(i),\ i\in \Z^d\}$ be i.i.d. random variables independent
of the $\{X_i,\ i\in \N\}$, with mean 0 and finite variance $\sigma^2$.
Kesten and Spitzer showed in~\cite{kesten-spitzer} that in dimension 1, the
following weak convergence in law (over both randomness) holds 
\[
\frac{1}{n^{3/4}}\sum_{k=1}^{[nt]} \xi(S_k)
\xrightarrow{law} \Delta_t,
\]
where $\Delta_t$ is a non-Gaussian, self-similar process of order $3/4$
with stationary increments. When $d=2$, Bolthausen~\cite{bolthausen}
established that
\[
\frac{1}{\sqrt{n\log(n)}} \sum_{k=1}^{[nt]} \xi(S_k)\xrightarrow{law}
\frac{\sigma}{\sqrt{\pi} {\tt det}(\Sigma)^{1/4}} B_t.
\]
When $d\ge 3$, Kesten and Spitzer essentially established
(in~\cite{kesten-spitzer}, page 10) that
\[
\frac{1}{\sqrt n} \sum_{k=1}^{[nt]} \xi(S_k)\xrightarrow{law}
\frac{\sigma}{\sqrt{E[N(0)]}} B_t,
\]
where $E[N(0)]$ is the expected number of visits to the origin
of the (transient) random walk $\{S_n,n\in \N\}$. 
Our interest was to understand how these super-diffusive scaling would
reflect in the large deviation speed rates.
The use of Brownian motion, rather than random walks, has
technical advantages: on one hand, we have the spectral analysis
and the classical estimates for Schr\"odinger semi-groups
at our disposal, and on the other hand, we have a clean scaling property.

Some Large Deviations estimates for $\frac{1}{t}\int_0^t \xi(B_s)ds$ 
were obtained in~\cite{remillard} for the annealed case. In particular
the speed rate was obtained in dimension 1, but not the rate functional.
Besides, in $d>1$, not even the correct speed was discovered.

We now give some heuristics to explain the correct speed rates
in estimating, for any real $y$, the probability of the event
$\AA \triangleq \acc{ (\xi,B)\, : \, \bra{L_t,\xi} \approx y}$
where $L_t = \frac{1}{t}\int_0^t \delta_{B_s} \, ds$ is
the occupation measure of Brownian motion, and $\bra{\cdot,\cdot}$
is the duality bracket between measures and functions.
By scale invariance, we have for any $r>0$, $B_{r^2s}
\overset{law}{=}rB_s$. Thus, we have to find the probability of 
the event 
\begin{equation} 
\label{eq0.2}
\AA
= \acc{ (\xi,B) \, : \, \, \bra{L_{\frac{t}{r^2}}, \xir}
\approx y } \, , 
\end{equation}
 where  
\begin{equation} 
\label{eq0.3}
\xir(x)\triangleq \xi(rx) = \sum_{i \in Z^d} \xi(i) 
\ind_{Q_i(\frac{1}{r})}(x), 
\mbox{ and } Q_i(\frac{1}{r}) \triangleq
\frac{i}{r}+[-\frac{1}{2r},\frac{1}{2r}[^d.
\end{equation}
Thus, the Brownian scale invariance
has allowed us to ``coarse-grain'' the field.
Indeed, we think of $\acc{\xi(i),i\in \Z^d}$ as our microscopic description
and introduce the empirical density $\xir$, which represents
coarse graining over about $r^d$ sites.

Now, a Large Deviations Principle (LDP)
holds for the field $\xir$ integrated against continuous functions
with compact support (see e.g.~\cite{bkl}). In other words, for any $y>0$, 
$\varphi \in \CC_c(\R^d)$, and $r$ large
\begin{equation} 
\label{eq0.6}
\frac{1}{r^d} \log \displaystyle{\otimes_{i\in \Z^d} \nu_i}
\cro{\bra{\xir, \varphi} \approx y}
\approx -\!\!\!
\inf_{u: |u(x)|\le 1} \acc{ I(u): \! \! \int_{\R^d}\!\!
u(x)\varphi(x)\, dx\approx y}.
\end{equation} 
with 
\begin{equation} 
\label{eq0.7}
I(u)\triangleq \int_{\R^d} H(u(x))dx,
\quad {\tt and} \quad
H(x)\triangleq \sup_{y}\acc{xy-\log E_{\nu}\cro{\exp(y\xi(0))}}.
\end{equation} 
On the other hand, the Donsker-Varadhan
theory provides a LDP for the occupation measure $L_{t/r^2}$
in the weak topology.

Thus, when we average with respect to both randomness, 
i.e. in the annealed case, 
it is natural to look for a LDP  by 
a contraction principle (cf.~\cite{dembo-zeitouni}).
Assume for a moment that we are entitled to do so. 
Then, the correct speed appears as
one equals $t/r^2$ with $r^d$, i.e. as one equals the speed rates
for each marginal LDP. Thus, this yields the correct speed $t^{d/(d+2)}$.
Moreover, the rate function is
\[
\II(y) = \inf \acc{I(u)+\LL(\mu)\, : \,\, \bra{\mu,u}=y} \, ,
\]
where $\LL$ is the rate function for the Brownian occupation measure. 
However, when using a contraction principle, we face two problems. 
First, the map $(u,\mu)\mapsto \int u d\mu$ is not continuous
in the product of the weak topologies. The remedy is to regularize the
field: if $\{\pd\}$ is an approximate identity, one first has to 
replace $\xir$ by $\pd*\xir$. Second, the LDP for the Brownian
occupation measure is a weak one, i.e. the upper bound
is only valid for compact sets. The standard
trick, which has first been used by Donsker and Varadhan \cite{sausage},
is to replace the Brownian motion
by a process for which we have a ``full'' LDP, for instance
the Brownian motion on the torus $\TT(A)$ of side $A$. This compactification
is possible in our situation, since we show that if we integrate
first with respect to the law of $\xir$,
\[\P\cro{\bra{L_{\frac{t}{r^2}},\pd*\xir} \geq y }
\leq e^{-r^d F_A(L^A_{t/r^2})} \, ,
\ {\tt with}\ 
F_A(\mu)=\inf \acc{ I_A(u): \, \bra{\mu,\pd*u} \geq y} \, ,
\]
where $L^A$ is the occupation measure for the Brownian on $\TT(A)$,
and $I_A$ has the same expression as $I$ in \refeq{eq0.7} with
$\TT(A)$ instead of $\R^d$. 
The upper bound follows then from Varadhan's integral lemma
and coincides with the lower bound.

Another standard way to obtain a LDP is to use G\"artner-Ellis method,
i.e. to look for the asymptotics of the log-Laplace transform 
\[
\frac{t}{r^2} \log \tilde{E}_0
\cro{ \exp(\frac{\alpha}{r^2} \int_0^t \xi(B_s) \, ds)}
=  \frac{t}{r^2} \log \tilde{E}_0
\cro{\exp(\alpha \int_0^{t/r^2} \xir(B_s) \, ds)} \, ,
\]
where $\tilde{E}_0$ denotes the annealed law.
By Feynman-Kac formula, this behavior is related to the 
(annealed) behavior of the principal eigenvalue of
the random operator $- \frac{1}{2} \triangle - \alpha \xir$.
Similar quantities have been thoroughly studied both in the
annealed and in the quenched setting, for different kinds of
potential $\xi$ and different scaling $r$: for instance
Sznitman \cite{sznitman}, Merkl \& W\"uthrich \cite{merkl-wuthrich-ptrf}
\cite{merkl-wuthrich-ihp} \cite{merkl-wuthrich-spa} 
for the case of a Poissonian potential; G\"artner \& Molchanov 
\cite{gartner-molchanov-cmp} \cite{gartner-molchanov-ptrf}, 
Biskup \& K\"onig \cite{biskup-konig-ap} for the i.i.d case;
G\"artner \& K\"onig \cite{gartner-konig-aap},  
G\"artner, K\"onig \& Molchanov \cite{gartner-konig-molchanov-ptrf} 
for more general potentials.  This method leads to a LD upper bound
which is necessarily convex, being defined as a Legendre transform.
However, the functional
$\II$ is not convex in general, and this method is doomed to fail.

What about the case with a fixed field $\xi$,
i.e. the quenched case? First, note that if $\sigma(R)$
denotes the Brownian exit time from a cube of radius $R$, then
by classical results, there is a constant $C$ such that
\begin{equation}
\label{eq0.8}
P_0 \cro{\sigma(R\frac{t}{r^2})\le \frac{t}{r^2}} \le C\exp \pare{- 
\frac{R^2}{2}\frac{t}{r^2}}.
\end{equation}
Hence, we can restrict everything to a box $Q$ of size $Rt/r^2$. 
Now, to establish estimates
holding $\xi$-almost surely, a pattern of the scaled field $\xir$
(on a macroscopic domain) should persist as we take $t$ to infinity.
By a Borel-Cantelli argument, this happens as soon  as $r^d=\log(t/r^2)$.
Indeed, the cost for $\pd*\xir$ to look like
a definite profile $u$ on a unit cube, is of order $\exp(-r^d I(u))$.
Since the smoothed empirical density
$\pd*\xir$ is almost independent on the different cubes of a partition
of $Q(Rt/r^2)$ into unit cubes, the probability that in one of
the element of the partition, $\pd*\xir$ is close to $u$ is of order
\begin{equation} 
\label{eq0.9}
1- \pare{1-\exp\pare{-r^dI(u)}}^{(\frac{Rt}{r^2})^d}.
\end{equation} 
This is almost 1, if $r^d=\log(t/r^2)$, whose root we call $r_t$,
and $I(u)<d$.
Now, forcing the Brownian motion to stay in a unit cube during a
time $t/r_t^2$ costs of the order of $\exp(-ct/r_t^2)$. Thus, we 
have the heuristic speed $t/r_t^2$ with $r_t^d\triangleq\log(t/r_t^2)$.
Following this strategy and optimizing over all admissible profiles
$u$, we obtain 
\[
\begin{array}{ll}
& P_0\cro{\bra{L_{t/r_t^2};\xir} \approx y}
\geq \exp\pare{-\frac{t}{r_t^2} \JJ_1(y)} \, , 
\\ 
\mbox{ where }
& \JJ_1(y) = \inf_{u,\mu} \acc{\LL(\mu)\, : \,
\bra{\mu,u}=y \, , \, I(u) < d} \, .
\end{array} 
\]

The quenched large deviations upper bound is obtained using G\"artner-Ellis
method. As already mentioned, we are led to study the almost sure behavior 
of the principal eigenvalue of the random  Schr\"odinger operator
$-\frac{1}{2} \bigtriangleup - \alpha \xir$ 
with boundary Dirichlet conditions on 
$Q(Rt/r_t^2)$. In the case of a Poissonian potential, this study has
been carried out by Merkl \& W\"uthrich \cite{merkl-wuthrich-ihp}.
We rely here  on  a localization lemma borrowed 
from G\"artner and K\"onig \cite{gartner-konig-aap}, which has also been
crucial in the papers \cite{biskup-konig-ap}, 
\cite{gartner-konig-molchanov-ptrf}, \cite{merkl-wuthrich-ihp}.  
According to this lemma,  the principal eigenvalue
is  close to the minimum of the principal eigenvalues of the
same operators over boxes of fixed size forming a partition of 
$Q(Rt/r_t^2)$.  
This leads to a quenched upper bound with a rate functional $J$ which
is convex.

The problem is now to identify $\JJ_1$ with $J$. It is easy
to check that $J$ is the greatest convex minorant of $\JJ_1$. 
However, the convexity of $\JJ_1$ could not be established. Hence, 
we use an approach developed 
in~\cite{asselah-castell}, and 
we convexify through
a sequence of scenarios: the $n$-th one corresponds to
partitioning $[0,T]$ into  $n$  time intervals, in each of which
the Brownian motion goes fast to a region
where the field $\pd*\xir$ has a
fixed deterministic profile, and stays there during this time interval.
To each scenario corresponds a lower bound of the type
\[
\lim_{\epsilon \rightarrow 0}
\varliminf_{t \rightarrow \infty}
\frac{r_t^2}{t} \log P_0 \cro{ \frac{r_t^2}{t}\int_0^{t/r_t^2}
\pd*\xir(B_s)ds\approx y}
\geq - \JJ_n(y) \, .
\]
The family of functions $\JJ_n$ is decreasing, and satisfies for any
$y_1,y_2$ and $\lambda\in ]0,1[$
\[
\lambda \JJ_n(y_1)+(1-\lambda)\JJ_n(y_2)\ge
\JJ_{2n}(\lambda y_1+(1-\lambda)y_2).
\]
Thus, the limit $\JJ(y) \triangleq \lim_{n \rightarrow \infty} \JJ_n(y)$
is convex. This enables us to identify $\JJ$ with the upper bound $J$.

Though we restrict ourselves to the i.i.d case, 
the crucial assumptions are that the rescaled
field $\xir$ is bounded and satisfies a LDP. 

The paper is organized as follows. In section \ref{cadre}, we
introduce the notations and state the main results. In
section \ref{anneal}, we prove the LDP for the annealed case.
In section \ref{quench}, we establish the LDP for the quenched case.
Section \ref{lemmes} gathers the proof of some technical lemmas.

\vspace{0,5cm}
\section{Notations and results.}
\label{cadre}

{\bf The random scenery.} Let $\{\xi(j), j \in \Z^d\}$ be
a family of i.i.d random variables  with values in $\R$.
We denote by $\P=\otimes_j\nu_j$ the law of the environment. Expectation
with respect to $\P$ is denoted by $\E$.  We assume that 
\begin{equation} 
\label{hyp-champ}
\P\mbox{-a.s. } -1 \leq \xi(0) \leq 1\, \quad 
\E\cro{\xi(0)} = 0 \, , \quad \mbox{ and } \quad
\E\cro{\xi(0)^2} \neq  0 \, . 
\end{equation} 
We will denote $m \triangleq \mbox{essinf}(\xi(0))$, and 
$M \triangleq \mbox{esssup}(\xi(0))$. 

Let $\Lambda$ be the log-Laplace transform of $\xi(0)$:
\begin{equation} 
\label{eq1.1}
\forall \alpha \in \R, \quad \Lambda(\alpha) \triangleq  
\log \E\cro{e^{\alpha \xi(0)}} \, . 
\end{equation}
$\Lambda$ is convex, everywhere finite by \refeq{hyp-champ}. 
Moreover, since $\E(\xi(0))=0$, $\Lambda(\alpha) \geq 0$, and $\Lambda(0)=0$.
Let $H$ be the Legendre transform of $\Lambda$:
\begin{equation} 
\label{eqi.2}
\forall y \in \R, \quad  H(y)\triangleq 
\sup_{\alpha \in \R} \pare{\alpha y - \Lambda(\alpha)} \, .
\end{equation}
$H$ is convex, takes positive values, is increasing 
on $\R^+$, decreasing on $\R^-$. $H(0)=-\inf \Lambda =0$. 
$H(y) = + \infty$ for $y \notin [m,M]$, and $H(y) < \infty$ for
$y \in ]m,M[$.

If $Q(A) \triangleq [-\frac{A}{2};
+\frac{A}{2}]^d$, let $\MM(Q(A))$ (resp. $\MM_1(Q(A))$, 
$\MM_1^0(Q(A))$) 
be the set of finite signed measures on $Q(A)$ (resp. the
set of probability measures on $Q(A)$, the set of probability
measures with compact support included in  $Q(A)$),  
endowed with the topology of weak convergence (i.e the topology defined
by duality against continuous and bounded test functions).
For all $r > 0$, let $\xir$ be the function 
defined by 
\begin{equation}
\label{xir-def}
\xir(x)  \triangleq  \xi([rx])\, \,\,
\mbox{ where } [x] \mbox{ is the integer part of } x \, .
\end{equation}
$\{\xir(x), x \in Q(A)\}$ are then random variables with values in 
\[
\BB_1(A)=\acc{u \in L_{\infty}(Q(A)), \nor{u}_{\infty} \leq 1} \, .
\]
$\BB_1(A)$ will  be viewed as the subspace of $\MM(Q(A))$ of measures
whose density with respect to Lebesgue measure belongs to $\BB_1(A)$. 

A key result is the following
large deviations principle (see for instance \cite{bkl}).
\begin{lemma}
\label{lem.champ-LDP}
For all $A >0$, let us define the rate function $I_A$ on $\BB_1(A)$ by
\begin{equation} 
\label{rf-champ}
I_A(u) = \int_{Q(A)} H(u(x)) \, dx
\end{equation}
$I_A$ is convex, lower semi-continuous.

When $r \rightarrow \infty$, $\xir$ satisfies a LDP on $\BB_1(A)$,
with good rate function $I_A$ and speed $r^d$; i.e. for all measurable
subset $F$ of $\BB_1(A)$, 
\begin{equation} 
\label{eqi.3}
-  \inf_{u \in \stackrel{\circ}{F}} I_A(u)
\leq \varliminf_{r \rightarrow \infty} \frac{1}{r^d} \log \P( \xir \in F)
\leq \varlimsup_{r \rightarrow \infty} \frac{1}{r^d} \log \P(\xir \in F)
\leq -\inf_{u \in \bar{F}} I_A(u) \, .
\end{equation} 
\end{lemma}
For $A=\infty$, $I_A$ and $\BB_1(A)$  will simply be denoted by $I$
and $\BB_1$.

\vspace{.5cm}
\noindent 
{\bf The Brownian motion in random scenery.} Let $\{B_t, t \in \R^+\}$ be
a $d$-dimensional Brownian motion, independent of the random field 
$\xi$.  $E_x$ denotes expectation under the 
Wiener measure starting
from $x$. For $t > 0$, let $L_t \triangleq \frac{1}{t} \int_0^t \delta_{B_s}
\, ds$ be the Brownian occupation measure.
>From Donsker-Varadhan theory,
$L_t$ satisfies a weak LDP in $\MM_1(\R^d)$, with speed $t$,
and rate function
\[
\LL(\mu) = \left\{ \begin{array}{ll}
\displaystyle{\frac{1}{2}}
 \int \nor{\nabla \sqrt{\frac{d\mu}{dx}}}^2 \, dx \, ,
& \mbox{ if } \mu \ll dx \, : \\
+\infty & \mbox{ otherwise }.
\end{array}
\right.
\]
When $\mu$ is a measure, and $u$ is a function, $\bra{\mu,u}\triangleq
\int u \, d\mu$.
Our main interest in this paper is large deviations estimates
for the random additive functional $\bra{L_t, \xi}$
under the ``quenched'' measure $P_0$, and the ``annealed'' one
$\tilde{P}_0 \triangleq \E (P_0)$. Before describing our results,
we need more notations.  

In all the sequel, when  $D$ is a  domain of $\R^d$, $\CC^{\infty}_c(D)$
is the space of infinitely differentiable functions with compact support
in $D$. $H^1_0(D)$ is the Sobolev space obtained by closure 
of   $\CC^{\infty}_c(D)$ under the norm
\[ \nor{f}^2 = \int_D f^2(x) \, dx + \int_D \nor{\nabla f(x)}^2 \, dx  \, .
\]
When $V: D \mapsto \R$ is a bounded measurable function,  
 we will write $\lambda(V,D)$ for the principal
eigenvalue of the operator $-1/2 \triangle - V$, with Dirichlet
boundary condition on  $D$.
\[
\lambda(V,D) \triangleq \inf \acc{\frac{1}{2} \int_{D} 
\nor{\nabla f}^2(x) \, dx - \int_D V(x) f^2(x) \, dx :
f \in H^1_0(D), \int_D f^2(x) \, dx =1 } \, .
\]
For $V \equiv 0$, and $D=Q(1)$, $\lambda(V,D)$ will be denoted by
$\lambda_1(d)$.

\vspace{.5cm}
\noindent
{\bf The annealed large deviations principle.}
For any $y \in \R$, let us define
\begin{equation} 
\label{rf-annealed}
\II(y)\triangleq \inf\acc{I(u)+\LL(\mu): \,  u \in \BB_1, \mu \in \MM^0_1(\R^d),
\bra{\mu,u} = y} \, .
\end{equation}
Let $\tilde{\II}$ be the greatest lower semi-continuous minorant of $\II$:
\begin{equation}
\label{rfa.lsc}
\tilde{\II}(y) \triangleq 
\lim_{\epsilon \rightarrow  0} \inf_{|z-y|<\epsilon} \II(z) \, .
\end{equation}
\begin{theo} 
\label{lda.th}
Assume \refeq{hyp-champ}.  Then,  for any measurable subset  
$F$ of $\R$, 
\begin{equation} 
\label{uba.eq} 
\limsup_{t \rightarrow \infty} \frac{1}{t^{\frac{d}{d+2}}}
\log \tilde{P}_0 \cro{\bra{L_t,\xi}  \in F}
\leq - \inf_{y \in \bar{F}} \, \tilde{\II}(y) \,\, , 
\end{equation} 
\begin{equation} 
\label{lba.eq}
\liminf_{t \rightarrow \infty}  \frac{1}{t^{\frac{d}{d+2}}}
\log \tilde{P}_0 \cro{\bra{L_t, \xi}  \in F}
\geq - \inf_{y  \in  \stackrel{\circ}{F}} 
	  \, \tilde{\II}(y) \,\, .
\end{equation} 
$\tilde{\II}: \R \mapsto \R^+$ is lower semi-continuous, 
increasing on $\R^+$, decreasing on $\R^-$.
$\tilde{\II}(0)=0$, $\tilde{\II}(y) < \infty$ 
for $y \in ]m;M[$, 
$\tilde{\II}(y) = \infty$ 
for $y \notin [m;M]$. Moreover, for $ d \leq 4$, 
\begin{equation}
\label{eq0-rfa}
\liminf_{y \rightarrow 0} \frac{\tilde{\II}(y)}{|y|^{\frac{4}{2+d}}} > 0 \, .
\end{equation} 
\end{theo}

{\bf Remark.} For $d \leq 4$, $\lim_{y \rightarrow 0} 
\tilde{\II}(y)/|y|^{4/(2+d)}\in ]0, +\infty[$. 
We will not prove this fact, since our interest
is to show that $\tilde{\II}$ is not convex in 
dimension $d=3$ and $d=4$, and $\refeq{eq0-rfa}$ is enough for that purpose.
Actually, $\frac{4}{2+d} < 1$ for $d=3,4$. Hence, 
if $\tilde{\II}$ is convex,   
$\tilde{\II}(y) =+\infty$ for any $y \neq 0$. This contradicts 
the fact that $\tilde{\II}$ is finite on $]m,M[$.
 
\vspace{.5cm}
\noindent 
{\bf The quenched large deviations principle.}
For any $\alpha \in \R$, let 
\begin{equation} 
\label{log-laplace-quenched}
l(\alpha) \triangleq
 \inf_{u \in \BB_1} \acc{\lambda(\alpha u,\R^d): I(u) \leq d} \, ,
\end{equation} 
and for any $y \in \R$
\begin{equation} 
\label{rf-quenched}
J(y) \triangleq \sup_{\alpha \in \R} \acc{ \alpha y + l(\alpha)} \, .
\end{equation} 
Then, 
\begin{theo} 
\label{ldq.th}
Assume \refeq{hyp-champ}. Let us define $r(t)$ by the relation
$t = r^2(t) \exp(r^d(t))$. Then, $\P$-a.s., for any measurable subset  
$F$ of $\R$, 
\begin{equation} 
\label{ubq.eq} 
\varlimsup_{t \rightarrow \infty} \frac{r^2(t)}{t}
\log P_0 \cro{\bra{L_t,\xi}  \in F}
\leq - \inf_{y \in \bar{F}} \, J(y) \,\, , 
\end{equation} 
\begin{equation} 
\label{lbq.eq}
\varliminf_{t \rightarrow \infty} \frac{r^2(t)}{t}
\log P_0 \cro{\bra{L_t, \xi}  \in F}
\geq - \inf_{y \in \stackrel{\circ}{F}} 
   \, J(y) \,\, .
\end{equation} 
$J: \R \mapsto \R^+$ is convex, lower semi-continuous, increasing on
$\R^+$, decreasing on $\R^-$. $J(0)=0$, $J(y) < \infty$ for any 
$y \in ]m,M[$, $J(y) = \infty$ for $y \notin [m;M]$. 
\end{theo}


\section{Annealed Bounds.}
\label{anneal}
In section~\ref{smooth}, we regularize the field. 
We prove the annealed LD lower bound in section \ref{alb}, and
the corresponding upper bound in section \ref{aub}. 
In all the sequel, we set for convenience $\tau=t/r^2$.


\subsection{Smoothing the field.}
\label{smooth}
Let $\psi$ be a rotationally invariant, nonnegative, 
smooth function with support in $Q(1)$
and integral 1. For $\delta>0$, let $\pd(x) \triangleq
 \psi(x/\delta)/\delta^d$. We denote by $*$  the convolution, that is
$u*v(x)=\int_{\R^d} u(x-y)v(y)dy$. 

\begin{lemma}
\label{reg-lem} 
For any $\epsilon >0$,
\begin{equation} 
\label{eq.lem1.1}
\lim_{\delta \to 0}\lim_{\tau\to\infty} \frac{1}{\tau} \log \sup_{u\in \BB_1}
P_0\pare{|\bra{L_{\tau},\pd*u-u}|>\epsilon }= -\infty ,
\end{equation} 
and 
\begin{equation} 
\label{eq.bis-lem1.1}
\lim_{\delta \to 0}\lim_{\tau\to\infty} \frac{1}{\tau} \log
\tilde P_0\pare{|\bra{L_{\tau},\pd*\xir-\xir}|>\epsilon }= -\infty.
\end{equation} 
\end{lemma}  

\vspace{.5cm}
\noindent
{\bf Proof.} 
In view of the classical fact
\begin{equation} 
\label{exit-time}
\lim_{R\to \infty}\lim_{\tau\to\infty}\frac{1}{\tau} \log P_0\pare{
\sigma(R\tau)< \tau}=-\infty,
\end{equation} 
the result (\ref{eq.lem1.1}) follows as soon as we show that
\[
\lim_{\delta \to 0}\lim_{\tau\to\infty} 
\frac{1}{\tau} \log \sup_{u\in \BB_1}
P_0\pare{\AA_{reg}(\tau,\delta, u)}=-\infty \, ,
\] 
where 
\begin{equation} 
\label{eq1.4}
\AA_{reg}(\tau,\delta ,u)\triangleq \{B_.:\,  
\va{\bra{L_{\tau},\pd*u-u}}>\epsilon ,
\sigma(R\tau)\ge \tau\}.
\end{equation} 
We  only estimate the probability of the event
$\AA(\tau,\delta,u) \triangleq 
\{B_.:\ \bra{L_{\tau},\pd*u-u}>\epsilon,
\sigma(R\tau)\ge \tau\}$, and the remaining
part of $\AA_{reg}(\tau,\delta ,u)$ can be dealt with similarly.

By Chebychev's inequality, we have for any $a>0$
\begin{equation} 
\label{eq1.5}
P_0(\AA(\tau,\delta,u))\leq E_0\cro{\exp\pare{a\int_0^{\tau}
\pare{\pd*u-u}(B_s)\, ds}\ind_{\sigma(R\tau)>\tau}}e^{-\tau a\epsilon}.
\end{equation} 
Using classical bounds (see e.g. Theorem 3.1.2, p.93 of~\cite{sznitman}),
there is $c(d)>0$ such that
\begin{equation} 
\label{eq1.5bis}
P_0(\AA(\tau,\delta,u))\le  c(d)
(1+(\tau \lambda(a(\pd*u-u),Q(R\tau)))^{d/2})
e^{-\tau\pare{a\epsilon +\lambda(a(\pd*u-u),Q(R\tau))}},
\end{equation} 
Note that when $u \in \BB_1$,
\[
\lambda(a(\pd*u-u),Q(R\tau))\le 2a +\lambda_1(d)/(R\tau)^2. 
\]
Moreover, we prove in section \ref{lemmes} the following result.

\begin{lemma}
\label{lem-tech3} For any $u\in \BB_1(\R^d)$, $\mu\in \MM_1(\R^d)$,
$\delta >0$ and $\epsilon_1>0$,  we have
\begin{equation} 
\label{eq.tech3}
|\int(\pd*u-u)d\mu|\le \frac{c_0\delta ^2}{\epsilon_1} \LL(\mu)+2\epsilon_1.
\end{equation} 
\end{lemma}

Thus, for any $a > 0$,
\begin{equation} 
\label{eq1.10}
P_0(\AA(\tau,\delta,u))
\leq \pare{1+ \pare{2a\tau+\frac{\lambda_1(d)}{R^2\tau}}^{d/2}}
e^{-\tau\pare{a\epsilon +\inf_{\mu\in \MM_1^0(Q(R\tau))}
\acc{\LL(\mu)(1-\frac{a c_0 \delta ^2}{\epsilon_1})-2a\epsilon_1}}}
\, .
\end{equation} 
We set $\epsilon_1=\epsilon/4$, and $a = \frac{\epsilon_1}{2 c_0 \delta^2}$. 
We have
\begin{equation}
\label{reg-ubq}
\sup_{u\in \BB_1} P_0(\AA(\tau,\delta,u))
\leq 
\pare{1+ \pare{\frac{\epsilon \tau}{4 c_0 \delta^2}
	+ \frac{\lambda_1(d)}{R^2\tau}}^{d/2}}
\exp\pare{-\tau \pare{\frac{\epsilon^2}{16 c_0 \delta^2}+
			\frac{\lambda_1(d)}{2R^2\tau^2}}}.
\end{equation} 
Hence, we also have
\begin{eqnarray} 
\label{reg-uba}
\lefteqn{
\tilde{P}_0 \cro{ \bra{L_{\tau}, \pd * \xir - \xir} \geq \epsilon \, :
\sigma(R \tau) > \tau} } & & \nonumber \\
&& \leq 
\pare{1+ \pare{\frac{\epsilon \tau}{4 c_0 \delta^2}
	+ \frac{\lambda_1(d)}{R^2\tau}}^{d/2}}
\exp\pare{-\tau \pare{\frac{\epsilon^2}{16 c_0 \delta^2}+
			\frac{\lambda_1(d)}{2R^2\tau^2}}}.
\end{eqnarray} 

Equations \refeq{eq.lem1.1} and \refeq{eq.bis-lem1.1}  follow.
\hspace*{\fill} \rule{2mm}{2mm}

\subsection{The annealed lower bound.}
\label{alb}
By Lemma~\ref{reg-lem}, we can now replace $\bra{L_{\tau}, \xir}$ by
$\bra{L_{\tau}, \pd * \xir}$. The aim of regularizing
is the following.
\begin{lemma}
\label{cont-bra}
For any  $A > 0$,  the function $(\mu, u)\in \MM^1(Q(A)) 
\times \BB_1(A) \mapsto \bra{\mu, \pd * u}$, is continuous
in the product of weak topologies.
\end{lemma}

The proof of this lemma is given in section \ref{lemmes}. 
Since $(L_{\tau},\xir)$ satisfies
a LDP in the product of weak topologies, we immediately get 
the LD lower bound:

\begin{lemma}
\label{alb-lem}
Let $r$ and $\tau$ be such that $\tau=r^d$. Then,  for any $y \in \R$ 
and any  $\epsilon > 0$,
\[
\lim_{\epsilon \rightarrow 0} \, 
\liminf_{\tau \rightarrow \infty}
\frac{1}{\tau} \log \tilde{P}_0 
\cro{ |\bra{L_{\tau},\xir} -y|<\epsilon}
\geq - \tilde{\II}(y) \, ,
\]
where $\tilde{\II}(y)$ is defined in \refeq{rfa.lsc}.
\end{lemma} 

\vspace{.5cm}
\noindent
{\bf Proof}: 
From
\[
\tilde{P}_0 \cro{|\bra{L_{\tau},\xir}-y|<\epsilon}
\geq \tilde{P}_0 \cro{\bra{L_{\tau},\pd * \xir} -y|<\frac{\epsilon}{2}}
- \tilde{P}_0 \cro{ \va{\bra{L_{\tau},\pd * \xir -\xir}} \geq 
\frac{\epsilon}{2}} \, ,
\]
and Lemma~\ref{reg-lem}, it is enough to prove that for any $y \in \R$,
\begin{equation}
\label{alb-reg}
\lim_{\epsilon \rightarrow 0} \liminf_{\delta \rightarrow 0} 
\liminf_{\tau \rightarrow \infty}
\frac{1}{\tau} \log \tilde{P}_0 \cro{|\bra{L_{\tau}, \pd * \xir}-y|<\epsilon}
\geq - \tilde{\II}(y) \, .
\end{equation}
Let $A, \epsilon, \delta$ be fixed positive numbers. 
Let $\mu_0 \in \MM^0_1(Q(A))$, and $u_0 \in \BB_1(A)$ be such that
$|\bra{\mu_0, \pd * u_0} -y|<\frac{\epsilon}{2}$. By Lemma~\ref{cont-bra}, one
can then find a weak neighborhood $\VV_1(\mu_0)$ of $\mu_0$, and
a weak neighborhood $\VV_2(u_0)$ of $u_0$, such that $\forall 
\mu \in \VV_1$, $\forall u \in \VV_2$, $|\bra{\mu, \pd * u}-y|<\epsilon$. 
Hence,
\begin{equation} 
\tilde{P}_0 
\cro{|\bra{L_{\tau}, \pd * \xir}-y|<\epsilon}  
\geq 
\tilde{P}_0
\cro{ L_{\tau} \in \VV_1; \xir \in \VV_2} 
=  P_0 \cro{ L_{\tau} \in \VV_1}
\P \cro{\xir \in \VV_2} \, ,
\end{equation}
by independence. Applying now the LDP for $L_{\tau}$ and $\xir$, we get
for $\tau = r^d$,
\begin{equation}
\liminf_{\tau \rightarrow \infty} \frac{1}{\tau} \log  \tilde{P}_0 
\cro{|\bra{L_{\tau}, \pd * \xir}-y|<\epsilon} 
\geq - \LL(\mu_0) - I_A(u_0) \, .
\end{equation}
Taking the supremum over admissible $(\mu_0, u_0)$ leads to
\begin{equation} 
 \liminf_{\tau \rightarrow \infty} \frac{1}{\tau} \log  \tilde{P}_0 
\cro{|\bra{L_{\tau}, \pd * \xir}-y|<\epsilon} 
\geq - \inf\acc{ \II_{A,\delta}(z): |z-y|<\frac{\epsilon}{2}} \, ,
\end{equation}
where
\[
\II_{A,\delta}(z) \triangleq 
\inf \acc{ \LL(\mu) + I_A(u): u \in \BB_1(A), \mu \in \MM^0_1(Q(A)),
\bra{\mu, \pd * u} =z}.
\]
Now, it is easy to see that 
\begin{eqnarray} 
& & \limsup_{A \rightarrow \infty} \inf_{|z-y|<\frac{\epsilon}{2}}
\II_{A,\delta}(z) \leq \inf_{|z-y|<\frac{\epsilon}{2}}
\II_{\delta}(z) \, , \nonumber \\
& \mbox{where} & 
\II_{\delta}(z) \triangleq  \inf \acc{\LL(\mu) + I(u):
\mu \in \MM^0_1(\R^d), u \in \BB_1, \bra{\mu,\pd * u} =z} \, .
\label{def-Id}
\end{eqnarray} 
We now take $\delta$ to 0. In view of Lemma~\ref{lem-tech3}, 
it is easy to see that
\begin{equation}
\limsup_{\delta \rightarrow 0} 
\inf_{|z-y|<\frac{\epsilon}{2}} \II_{\delta}(z) 
\leq \inf_{|z-y|<\frac{\epsilon}{4}} \II(z) \, .
\end{equation}
This ends the proof of \refeq{alb-reg} and of Lemma~\ref{alb-lem}.
\hspace*{\fill} \rule{2mm}{2mm}

\subsection{The annealed upper bound.}
\label{aub}
We now prove the annealed upper bound. Note that
by assumption $\nor{\xir}_{\infty} \leq 1$, so that
$\va{\bra{L_{\tau},\xir}} \leq 1$. Hence, it is enough to prove 
the weak large deviations upper bound. By Lemma~\ref{smooth}, the
problem is thus  reduced to prove that
\begin{equation} 
\label{alb.eq}
\lim_{\epsilon \rightarrow 0}  \, 
\limsup_{\delta \rightarrow 0} \, 
\limsup_{\tau \rightarrow \infty}
\frac{1}{\tau}
\log \tilde{P}_0 
\cro{|\bra{L_{\tau},\pd*\xir}-y|<\epsilon}
\leq - \tilde{\II}(y) \, .
\end{equation} 
In contrast with the lower bound, we cannot obtain \refeq{alb.eq}
by contraction, since $L_{\tau}$ does not satisfy a full LDP.
We begin by the following lemma.

\begin{lemma} 
\label{+alb}
Let $r$ and $\tau$ be such that
$\tau=r^d$. Then for any $y >0$, 
for any $\epsilon > 0$,
\[
\limsup_{\delta \rightarrow 0} \,
\limsup_{\tau \rightarrow \infty}
\frac{1}{\tau} \log \tilde{P}_0
\cro{ \bra{L_{\tau},\pd * \xir} \geq y}
\leq - \inf_{z \geq y-2 \epsilon} \II(z) \, .
\]
\end{lemma} 

\vspace{.5cm}
\noindent
{\bf Proof.} Let $y, \epsilon > 0$ be fixed.
$P_0$-a.s, for any $a > 0$,
\begin{eqnarray*}
\P\cro{ \bra{L_{\tau},\pd * \xir} \geq y}
&\leq& \exp(-a \tau y ) 
\E\cro{\exp(a \tau \bra{\pd * L_{\tau}, \xir})}
\\
&= & \exp(-a \tau y)
\E\cro{\exp \pare{ a \tau \sum_{i \in \Z^d} \xi(i)
\int_{Q_i(\frac{1}{r})} \pd*L_{\tau}(x) \, dx}} 
\\
&= & \exp(-a \tau y) \exp\pare{ \sum_{i \in \Z^d} 
\Lambda(a \tau \int_{Q_i(\frac{1}{r})} \pd*L_{\tau}(x) \, dx)} \, .
\end{eqnarray*}
Since $\tau=r^d=1/|Q_i(\frac{1}{r})|$, we have by convexity of $\Lambda$
\[
\Lambda(a \tau \int_{Q_i(\frac{1}{r})} \pd*L_{\tau}(x) \, dx)
\leq \tau \int_{Q_i(\frac{1}{r})}    \Lambda(a \pd*L_{\tau}(x)) \, dx \, .
\]
Hence, $P_0$-a.s., $\forall a > 0$,
\[
\P \cro{ \bra{L_{\tau},\pd * \xir} \geq y}
\leq 
\exp(-a \tau y) 
\exp \pare{ \tau \int_{\R^d} \Lambda(a \pd*L_{\tau}(x)) \, dx} \, .
\]
Let $A$ be a large number which will be sent to infinity later.
We cover  $\R^d$ by boxes $\{Q_i(A), i\in \Z^d\}$ of diameter $A$.
The center of $Q_i(A)$ is $iA$. We get 
\[
\int_{\R^d} \Lambda(a \pd*L_{\tau}(x)) \, dx
= \int_{Q_0(A)} \sum_{i \in \Z^d}\Lambda(a \pd*L_{\tau}(x+iA)) \, dx \,  .
\]
Now, by H\"older's inequality, we have for any $x, y \geq 0$,
$\Lambda(x) + \Lambda(y) \leq \Lambda(x+y)$. Thus,
\[
\int_{\R^d} \Lambda(a \pd*L_{\tau}(x)) \, dx
\leq \int_{Q_0(A)} \Lambda(a \sum_{i \in \Z^d} \pd*L_{\tau}(x+iA)) \, dx \, .
\]  
$\pd$ being rotationally invariant, 
$\sum_i  \pd*L_{\tau}(x+iA) = \frac{1}{\tau} 
 \int_0^{\tau} \sum_i \pd(\va{x+iA-B_s}) ds$. For
$\delta < A/2$, there is at most one non-vanishing term 
among $\{\pd(|x+iA-B_s|),\ i\in \Z^d\}$, whose index $i$
is such that $|x+iA-B_s|=d_A(x_A,B^A_s)$, 
where $d_A$ denotes the Riemanian metric on the
torus $\TT(A)$ of diameter $A$,  $x_A$ and $B^A_s$ being the
projection of $x$ and $B_s$ on $\TT(A)$. Setting $\psi^A_{\delta}:
(x_A,y_A) \in \TT(A) \times \TT(A) \rightarrow \pd(d_A(x_A,y_A))$, 
we have thus proved that $P_0$-a.s., for any positive $a$, 
$A$, and $\delta < A/2$, 
\[
\P \cro{ \bra{L_{\tau},\pd * \xir} \geq y} 
\leq \exp(-a \tau y)  
\exp ( \tau \int_{\TT(A)} \Lambda(a \psi^A_{\delta} * L^A_{\tau}(x))
\, dx ) \, . 
\] 
Taking now the infimum in $a > 0$ yields $P_0$-a.s,  
$\forall A > 0$, $\forall \delta < \frac{A}{2}$,  
\[ 
\P \cro{ \bra{L_{\tau},\pd * \xir} \geq y} 
\leq \exp\pare{ - \tau F_A^{\delta}(y;L^A_{\tau})} \, ,  
\]   
where 
\[
\forall \mu \in \MM_1(\TT(A)),\qquad   
F_A^{\delta}(y;\mu) \triangleq \sup_{a > 0}\acc{ a y  
-\int_{\TT(A)} \Lambda(a \psi^A_{\delta} * \mu(x)) \, dx}.
\]
Note that $\mu \mapsto F_A^{\delta}(y;\mu)$ 
is nonnegative and lower semi-continuous
as the supremum of continuous functions. 
Moreover, $L^A_{\tau}$ satisfies
a full LDP with the good rate function $\LL_A$. Hence, integrating with 
respect to $E_0$, and applying Varadhan's integral lemma (see for instance 
Lemma~4.3.6 in \cite{dembo-zeitouni}), yields
\[
\limsup_{\tau \rightarrow \infty} 
\frac{1}{\tau} \log 
\tilde{P}_0 \cro{ \bra{L_{\tau},\pd * \xir}\geq y}  
\leq - \inf_{\mu \in \MM_1(\TT(A))} 
\acc{ \LL_A(\mu) + F_A^{\delta}(y;\mu) }
\, .
\]

We apply now the following lemma, which is proved in section \ref{lemmes}. 
\begin{lemma}  
\label{pb-dual} 
For any $A> 0$, $y >0$, and $f$ a continuous  probability density 
on $\TT(A)$, 
\[ 
\sup_{a > 0} \acc{ a y - \int_{\TT(A)} \Lambda(a f(x)) \, dx}
= \inf_{u \in \BB_1(A)} 
\acc{I_A(u): \int_{Q(A)} f(x) u(x) \, dx \geq y} \, .
\]
\end{lemma}
Thus, for any $A > 0$, and any $\delta < A/2$, 
\[
\limsup_{\tau \rightarrow \infty} 
\frac{1}{\tau} \log 
\tilde{P}_0 \cro{ \bra{L_{\tau},\pd * \xir} \geq y}  
\leq - \inf_{z \geq y}
\bar{\II}_{\delta,A}(z) \, ,
\] 
where we have defined  
\[ 
\bar{\II}_{\delta,A}(z) \triangleq \inf
\acc{ \LL_A(\mu) + I_A(u): \mu \in \MM_1(\TT(A)), u \in \BB_1(A), 
\int_{Q(A)} u \, \, d(\pd^A*\mu)  = z}
\, .
\]  

We consider now the limit $A$ to infinity. 

\begin{lemma} 
\label{lem2.2} Let $\II_{\delta}$ be defined by \refeq{def-Id}.
For any positive $y$ and $\epsilon$, 
\begin{equation}
\label{eq5.4}
\inf_{z \geq y- \epsilon } \II_{\delta}(z)
\leq
\limsup_{A \rightarrow \infty}  \, 
\inf_{z \geq y} \bar{\II}_{\delta,A}(z) \, .
\end{equation}
\end{lemma} 

The proof of Lemma~\ref{lem2.2} is given in section \ref{lemmes}.

We now take $\delta$ to 0. It follows from Lemma~\ref{lem-tech3} that 
\begin{equation} 
\label{eq5.10}
\inf_{z \geq y- 2\epsilon } \II(z)
\leq
\liminf_{\delta \rightarrow  0} 
\inf_{z \geq y- \epsilon} \II_{\delta}(z) \, .
\end{equation} 
\hspace*{\fill} \rule{2mm}{2mm} 

We come now to the properties of the rate functional.
\begin{lemma}
Let $\tau$ and $r$ be such that $\tau=r^d$. Then, for any $y \in \R$, 
\begin{equation} 
\label{waub}
\lim_{\epsilon \rightarrow 0} \, \limsup_{\tau \rightarrow \infty}
\frac{1}{\tau} \log \tilde{P}_0 \cro{|\bra{L_{\tau},\xir}-y|<\epsilon}
\leq - \tilde{\II}(y) \, .
\end{equation} 
Moreover, 
$\tilde{\II}$ satisfies the properties listed in Theorem~\ref{lda.th}.
\end{lemma}

{\bf Proof.} Let us first prove the properties of
$\tilde{\II}$. Since $\LL$ and $I$  take positive values,
the same holds for $\II$ and $\tilde{\II}$. Also,
$\tilde{\II}$ is lower semi-continuous by
definition. Since  $H(0)=0$, taking $u \equiv 0$ 
in the infimum defining
$\II$ (see \refeq{rf-annealed}), 
we get $\II(0) \leq I(0)+ \inf \acc{ \LL(\mu): \mu \in \MM^0_1(\R^d)}=0$.
For any $y \in ]m,M[$, taking $u\triangleq y \ind_{Q(1)}$, and $\mu \in 
M^0_1(Q(1))$ in the infimum defining $\II$, leads to 
$\II(y) \leq H(y) + \lambda_1(d) < \infty$. Now if we assume that $
\II(y)$ is finite, one can find $\mu$ and $u$ such that 
$\bra{\mu,u}=y$, $\LL(\mu) < \infty$
and $\II(u) < \infty$. Hence $dx$-a.e., $u(x) \in [m, M]$. $\LL(\mu)$ being
finite, $\mu \ll dx$, and one therefore gets $y = \bra{\mu,u} \in [m,M]$. 

We prove \refeq{eq0-rfa} when $d\le 4$. For any $A>0$, 
we perform the change of variables,
\[ 
u_A(x) = u(Ax) \, , \quad d \mu_A = A^d \frac{d\mu}{dx}(Ax) dx \, ,
\]
and obtain
\begin{eqnarray} 
\II(y) &\!\! =\!\! & \inf_{u,\mu} 
\inf_{A > 0} \acc{ A^d I(u) + A^{-2} \LL(\mu):\ 
u \in \BB_1, \mu \in \MM^0_1(\R^d),\  \bra{\mu,u}=y} 
\nonumber \\
&\!\! =\!\! &  \pare{\frac{2}{d}}^{\frac{d}{d+2}} \pare{1 + \frac{d}{2}} 
\inf_{u,\mu}
\acc{ I(u)^{\frac{2}{d+2}} \LL(\mu)^{\frac{d}{d+2}}:\ 
u \in \BB_1, \mu \in \MM^0_1(\R^d),\ \bra{\mu,u}=y}  \, .
\nonumber \end{eqnarray} 
Now, there is $C > 0$ such that for  any
$x \in \R$, $H(x) \geq C x^2$. Hence, for another constant $C$, 
\begin{eqnarray*} 
\II(y) & \geq & C \inf_{u,\mu}
\acc{ \pare{\int y^2 u^2(x) \, dx}^{\frac{2}{d+2}} 
\LL(\mu)^{\frac{d}{d+2}}: \ u \in \BB_1, \mu \in \MM^0_1(\R^d),\ 
\bra{\mu,u}=1} 
\\
& = & C y^{\frac{4}{d+2}} \inf_{u,\mu} 
\acc{ \nor{u}_2^{\frac{4}{d+2}} \LL(\mu)^{\frac{d}{d+2}}:\ 
u \in \BB_1, \mu \in \MM^0_1(\R^d),\  \bra{\mu,u}=1}
\, .
\end{eqnarray*} 
It remains to prove that the infimum is strictly positive.
When $d \leq 4$, 
and for any $\mu \in \MM^0_1(\R^d)$ such that $\LL(\mu) < \infty$, we have
by a Nash type inequality (see for instance lemma 5 in 
\cite{castell-pradeilles}) that 
$\nor{\frac{d\mu}{dx}}_2 \leq C \LL(\mu)^{\frac{d}{4}}$. Hence,
for any $\mu,u$ such that $\bra{\mu,u} =1$, 
\begin{eqnarray*}
1=\bra{\mu,u} & \leq & \nor{u}_2 \nor{\frac{d\mu}{dx}}_2 \\
& \leq & C \nor{u}_2 \LL(\mu)^{\frac{d}{4}} \, .
\end{eqnarray*}  
This yields \refeq{eq0-rfa}.  

 Let us now prove the monotonicity of $\II$ (and thus of
 $\tilde{\II}$) on $\R^+$. Let $0 <  y_1 \leq y_2$. We can assume
 that $\II(y_2) < \infty$. Let then $\eta >0$, $\mu_2 \in \MM^0_1(\R^d)$,
 $u_2 \in \BB_1$ be such that $\bra{\mu_2,u_2}=y_2$ and 
 $\LL(\mu_2) + I(u_2) \leq \II(y_2) + \eta$. Let us define
 $\alpha \triangleq  y_1/y_2 
 \in ]0,1]$, $\mu_1 \triangleq \mu_2$, $u_1 \triangleq \alpha u_2$. Then
 $u_1 \in \BB_1$, and $\bra{\mu_1,u_1}=\alpha y_2 = y_1$. Hence,
 $\II(y_1) \leq \LL(\mu_1) + I(u_1) \leq \LL(\mu_2) + \alpha I(u_2)$,
 by convexity of $I$. Therefore, $\II(y_1) \leq \II(y_2) + \eta$
 for any  $\eta > 0$.

We now turn to the proof of \refeq{waub} for $y > 0$; the
negative case can be treated similarly. Choose $\epsilon > 0$ such that
$ 3 \epsilon < y$, then 
 \begin{eqnarray*} 
\varlimsup_{\tau  \rightarrow \infty}
\frac{1}{\tau} \log  \tilde{P}_0 \cro{|\bra{L_{\tau},\xir}-y|<\epsilon}
& &\leq \varlimsup_{\tau  \rightarrow \infty}
\frac{1}{\tau} \log  \tilde{P}_0 \cro{ \bra{L_{\tau},\xir} \geq y - \epsilon}
\\ 
 & &\leq - \inf_{z \geq y - 3 \epsilon } \II(z) \,
\mbox{  by Lemma~\ref{+alb}} \, ,
\\  
&  &\leq -  \inf_{|z-y|<3\epsilon} \II(z) \,
\mbox{  , since } \II \mbox{ is increasing on } \R^+ \, .
\end{eqnarray*}  
\hspace*{\fill}  \rule{2mm}{2mm}


\section{Quenched bounds} 
\label{quench}


\subsection{Quenched upper bound.}
The task at hand in this section is to prove \refeq{ubq.eq} of
Theorem~\ref{ldq.th}. Note that by assumption   \refeq{hyp-champ}, $\P$.a.s, 
$\forall t$, $\va{\bra{L_t,\xi}} \leq 1$. Therefore it is enough to 
prove the weak large deviations upper bound (i.e the upper bound for
compact sets). Using regularization
of the field (lemma \ref{reg-lem}), Brownian scaling, equation 
\refeq{exit-time} and the G\"artner-Ellis method, the problem is
reduced to study the large time asymptotics of 
\begin{equation}
\Lambda_{\tau }\pare{\alpha \pd  *\xir,Q(R \tau)}
\triangleq  E_0 \cro{
\exp \pare{ \int_0^{\tau} \alpha \pd  * \xir(B_s) \, ds};\,
\sigma_{R \tau} > \tau} \,\, ,
\end{equation} 
where as before $\tau=t/r^2$.
The next lemma gives the asymptotical behavior along some
subsequences.

\begin{lemma}
\label{limsup.Lambda.lem}
Let $\beta >1$, and let $(\tau_n)$ and $(r_n)$ be 
defined by $\tau_n= \exp(r_n^d) = \beta^n$.
Then $\forall \delta  > 0$, $\forall R >0$, 
$\forall \alpha \in \R$, $\P$-a.s.,  
\begin{equation}
\label{limsup.Lambda.eq}
\varlimsup_{n \rightarrow \infty} \frac{1}{\tau_n} 
\log \Lambda_{\tau_n}(\alpha \pd*\xin;Q(R \tau_n)) 
\leq - l(\alpha) \, ,
\end{equation}
where $l(\alpha)$ is defined by \refeq{log-laplace-quenched}.
\end{lemma} 

Before entering the proof, we note that Lemma~\ref{limsup.Lambda.lem}
is enough to prove the weak large deviations upper bound.
Indeed, it follows from the continuity of $l$ (see Lemma~\ref{prop.l})
and the bound $\nor{\pd *\xir}_{\infty} \leq 1$, that
we can make the ``$\P$-a.s.'' in the preceding lemma independent of
$\alpha$. By standard arguments (see for instance Theorem~4.5.3 of
\cite{dembo-zeitouni}), we have the weak upper bound 
along subsequences of the form 
$t_n^\beta = \log(\beta)^{2/d} n^{2/d} \beta^n$ with $\beta >1$.
Thus, for any $\beta > 1$,  we have $\P$-a.s., that for any $K$ compact,
\[
\limsup_{n \rightarrow  \infty} 
\frac{(r^{\beta}_n)^2}{t_n^{\beta}}
 \log P_0 \cro{\bra{L_{t^{\beta}_n},\xi} \in K}
\leq - \inf_{y \in K} J(y) \, ,
\]
The result follows now from the fact that for 
$t \in [t^{\beta}_n; t^{\beta}_{n+1}[$,
 \begin{equation}
\label{rest.suite}
 \va{\bra{L_t,\xi}-\bra{L_{t^{\beta}_n},\xi}}
 \leq 2 \frac{t^{\beta}_{n+1}-t^{\beta}_n}{t^{\beta}_n}
 \leq 2(\beta-1) + 2\beta \frac{C(d)}{n} \, ,
\end{equation} 
and the lower semi-continuity of $J$.

\vspace{.5cm}
\noindent  
{\bf Proof of Lemma~\ref{limsup.Lambda.lem}.}
In all the sequel, $R$, $\beta$, $\delta$ and $\alpha$ are fixed. 
For convenience, we will often not
mention the dependence in $n$ of $\tau$ and $r$.

\vspace{.5cm}
\noindent 
{\bf Step 1}. We begin to reduce the problem on boxes of fixed size,
using Lemma~4.6 of \cite{biskup-konig-ap}.  
First of all, note that as  in the proof of Lemma~\ref{reg-lem}, $\P$-a.s, 
\[
\Lambda_{\tau} \pare{\alpha \pd *\xir ,Q(R\tau)} 
 \leq C 
\pare{
1 + \pare{\tau 
\lambda (\alpha \pd  *\xir,Q(R\tau))}^{d/2}}
e^{- \tau 
\lambda (\alpha \pd  *\xir,Q(R\tau))} \, .
\]
Now, we cover $Q(R\tau)$ with $\NN \triangleq (1+[\frac{R\tau}{A}])^d$ boxes
of diameter $A$ whose centers we denote by $x_i,\ i=1,\dots,\NN$.
By Lemma~4.6 of \cite{biskup-konig-ap}, there is
constant $C$ such that for any $\alpha, R, r, \tau, A$ with $R \tau \geq A$,
\[
\lambda (\alpha \pd  *\xir,Q(R\tau))
\geq \min_{1 \leq k \leq \NN}
\lambda(\alpha \pd  *\xir, Q_k(A)) - \frac{C}{A^2}\, ,
\]
where we have used the  notation $Q_k(A)\triangleq x_k+Q(A)$. 
Since 
\[
\lambda (\alpha \pd  *\xir,Q(R\tau)) \leq
\frac{\lambda_1(d)}{(R\tau)^2}+|\alpha|,
\]
we are led to   
\begin{equation} 
\label{limsup.Lamba.minvp.eq}
\varlimsup_{\tau \rightarrow \infty} \frac{1}{\tau} \log 
\Lambda_{\tau} \pare{\alpha \pd  *\xir,Q(R\tau)} 
\leq 
\frac{C}{A^2}
- \varliminf_{\tau \rightarrow \infty}
\min_{1 \leq k \leq \NN} \acc{\lambda \pare{\alpha \pd  *\xir,
 Q_k(A)} } \, .
\end{equation}

\vspace{.5cm}
\noindent 
{\bf Step 2.}
We have now to estimate the $\P$-a.s. behavior of the minimum
in the above expression. The proof of the following lemma is
done in section \ref{lemmes}.

\begin{lemma} 
\label{fr.minvp.lem}
Let $r^d =\log(\tau)$ with $\tau=\beta^n$.  Then, for any
positive $\delta,A,R$, and any reals $\alpha,x$, we have
\begin{equation}
\label{fr.minvp.eq}
\begin{array}{l}
\varlimsup_{n \rightarrow \infty} \frac{1}{\log(\tau)}
\log \P \cro{
\min_{1 \leq k \leq \NN} 
\acc{\lambda \pare{\alpha \pd  *\xir,
 Q_k(A)}} \leq x} 
\\[.2cm]
\hspace*{2cm}
\leq d - \inf_{u \in \BB_1(A+1)}
\acc{ I_{A+1}(u): \, \lambda(\alpha \pd * u, Q(A+1)) \leq x}
\, .
\end{array}
\end{equation}
\end{lemma} 

Hence, for any $\epsilon > 0$ and $n$ sufficiently large, 
\begin{eqnarray*} 
\lefteqn{ \P \cro{ 
\min_{1 \leq k \leq \NN} \lambda(\alpha \pd  *\xir, Q_k(A)) \leq x}} \\
\\
& & \leq 
\exp\pare{
\log(\tau)\pare{d - \inf_{u \in \BB_1(A+1)}
\acc{ I_{A+1}(u): \, \lambda(\alpha \pd * u, Q(A+1)) \leq x} +\epsilon}}.
\end{eqnarray*} 
It follows from Borel-Cantelli lemma that
for any fixed $\alpha,\delta , A, R$, 
\begin{eqnarray} 
& & d < \inf_{u \in \BB_1(A+1)}
\acc{ I_{A+1}(u): \, \lambda(\alpha \pd * u, Q(A+1)) \leq x}
\nonumber 
\\
& &  \Rightarrow \label{liminf.vp.x}
\varliminf_{n \rightarrow \infty} \min_{1 \leq k \leq \NN}
\acc{\lambda \pare{\alpha \pd  *\xir, Q_k(A)}} > x \, ,
\, \P-\mbox{a.s.} \, .
\end{eqnarray}
We now prove that
\begin{eqnarray} 
& & x <  \inf_{u \in \BB_1(A+1)}
\acc{ \lambda(\alpha \pd*u, Q(A+1)): I_{A+1}(u) \leq d} 
\label{invbis.x} \\
& & \Rightarrow \label{inv.x}
\inf_{u \in \BB_1(A+1)} \acc{ I_{A+1}(u):
\lambda(\alpha \pd*u, Q(A+1))  \leq x} > d \, .
\end{eqnarray} 
Indeed, let $x$ satisfy (\ref{invbis.x}). For all
$u$ such that $\lambda(\alpha \pd*u, Q(A+1)) \leq x$, we have
then $I_{A+1}(u) > d$; in other words,
\begin{equation}
\label{inf.Ia}
\inf_{u \in \BB_1(A+1)} \acc{ I_{A+1}(u):
\lambda(\alpha \pd*u, Q(A+1))  \leq x} \geq d \, ,
\end{equation} 
with strict inequality if the infimum is reached, which is actually the case.
Indeed, it is proved in Lemma~\ref{cor-tech1} that 
$\acc{u \in \BB_1(A+1);  
\lambda(\alpha \pd *u, Q(A+1)) \leq x}$ is compact in weak topology.
Since  $I_{A+1}$ is lower semi-continuous in weak topology, $I_{A+1}$ reaches 
its minimum value on any compact set. 

 From \refeq{liminf.vp.x} and \refeq{inv.x}, we get for any
$\alpha, A, R$, that $\P$-a.s,
\begin{equation}
\label{liminf.vp.eq}
\varliminf_{n \rightarrow \infty} 
\min_{1 \leq k \leq \NN}
\acc{\lambda \pare{\alpha \pd  *\xir, Q_k(A)}}
\geq \inf_{u \in \BB_1(A+1)} 
\acc{\lambda(\alpha \pd * u, Q(A+1)): I_{A+1}(u) \leq d}
\,  \, .
\end{equation}

\vspace{.5cm}
\noindent
{\bf Step 3.} 
We show that for any $A > 0$,
\begin{equation}
\label{inf.IA.gt.inf.I.eq}
\inf_{u \in \BB_1(A)}
\acc{ \lambda(\alpha \pd*u, Q(A)): I_{A}(u) \leq d}
\geq 
\inf_{u \in \BB_1}
\acc{ \lambda(\alpha \pd*u, \R^d): I(u) \leq d} \, .
\end{equation}
Let $u_0 \in \BB_1(A)$ be such that $I_{A}(u_0) \leq d$. 
Let us consider the function $\tilde u_0 \in \BB_1$ defined by 
$\tilde u_0 \triangleq u_0 \ind_{Q(A)}$. Since $H(0) = 0$ , we have 
$I(\tilde u_0) = I_{A}(u_0) \leq d$. Therefore,
\[
\inf_{u \in \BB_1}
\acc{ \lambda(\alpha \pd*u, \R^d): I(u) \leq d} 
\leq \lambda(\alpha \pd*\tilde u_0, \R^d)
\leq \lambda(\alpha \pd*u_0, Q(A)) \, .
\]

Let us now prove that $\forall \delta   > 0$,
\begin{equation}
\label{inf.Id.gt.l}
\inf_{u \in \BB_1}
\acc{ \lambda(\alpha \pd*u, \R^d): I(u) \leq d}
\geq  l(\alpha) \, .
\end{equation} 
Indeed, by convexity of $H$, 
\[
I( \pd * u) = \int H(\pd * u) \leq 
\int \pd*H(u) = \int H(u) = I(u) \, .
\]
Therefore,
\[
\inf_{u \in \BB_1}
\acc{ \lambda(\alpha \pd*u, \R^d): I(u) \leq d} 
\geq \inf_{u \in \BB_1}
\acc{ \lambda(\alpha \pd*u, \R^d): I(\pd*u) \leq d}
\geq l(\alpha) \, .
\]

\vspace{.5cm}
\noindent
{\bf Step 4.} 
Lemma \ref{limsup.Lambda.lem} is then proved 
by putting \refeq{limsup.Lamba.minvp.eq}, 
\refeq{liminf.vp.eq}, \refeq{inf.IA.gt.inf.I.eq} and \refeq{inf.Id.gt.l}
 together, and by letting $A$ tend to infinity along
subsequences in \refeq{limsup.Lamba.minvp.eq}.
\hspace*{\fill} \rule{2mm}{2mm}


\subsection{Quenched lower bound.}
In this section, we prove \refeq{lbq.eq} of Theorem~\ref{ldq.th}. 

\vspace{.5cm}
\noindent 
{\bf Step 1. Almost sure  behavior of the field.}
\begin{lemma}
\label{a.s.comp.champ.lem}
Let $A > 0$ be fixed, and let $u \in \BB_1(A)$ be such that
$I_A(u) < d$. Let $\beta  > 1$ and let us define 
$\tau_n$ and $r_n $ by 
$\tau_n=e^{r_n^d}=\beta^n$. 
Then, for any positive $\delta$ and $\epsilon$, we have $\P$-a.s., 
that for $n$ sufficiently large, there is
a box $Q_k(\frac{[A r_n]}{r_n}) \subset 
Q(\tau_n/\log(\tau_n))$ such that 
\[
\nor{ \pd * \bar{\xi}_{r_n}
 - \pd * u_k}_{\infty, Q_k(\frac{[A r_n]}{r_n})}
< \epsilon      \, ,
\]
where $u_k$ denotes the translation of $u$ in the box
 $Q_k(\frac{[Ar_n]}{r_n})$.
\end{lemma}
{\bf Proof.} Let us note $A_r \triangleq \frac{[Ar]}{r} \approx A$ 
for large $r$. Define
\[
K \triangleq \acc{k \in \Z^d: 
Q_k(A_r) \subset Q(\tau/\log(\tau)) } \, , 
\] 
and let $\tilde{K}$ be the subset of $K$ corresponding to multi-integers
with even coordinates.
Note that as soon as $\delta + 1/r < A_r$, the functions
$\{\pd*\xir|_{Q_k(A_r)}; k \in \tilde{K}\}$ are independent. Moreover,
$A_r$ being an integer multiple of $1/r$, they also have the same
law. Therefore,
\[
\begin{array}{l}
\P \cro{\forall k \in K,  
\nor{ \pd * \xir - \pd * u_k}_{\infty, Q_k(A_r)}
\geq \epsilon  }
\\ 
\hspace*{1cm}  \leq \P \cro{\forall k \in \tilde{K},  
\nor{ \pd * \xir - \pd * u_k}_{\infty, Q_k(A_r)}
\geq  \epsilon }
\\
\hspace*{1cm}
 \leq  \P \cro{ 
\nor{ \pd * \xir - \pd * u }_{\infty, Q(A_r)}
\geq  \epsilon  }^{|\tilde{K}|} 
\\
\hspace*{1cm}
\leq \P \cro{ 
\nor{ \pd * \xir - \pd * u }_{\infty, Q(A)}
\geq  \epsilon  }^{|\tilde{K}|} \, .
\end{array}
\]
Now, $\acc{v \in \BB_1(A);  \nor{ \pd * v - \pd * u }_{\infty, Q(A)} < 
\epsilon }$
is an open neighborhood of $u$ in weak topology. Thus, by the LDP of
$\xir$,
\[ 
\varliminf_{r \rightarrow \infty}  \, 
\frac{1}{r^d} \log 
\P \cro{\nor{ \pd * \xir - \pd * u }_{\infty, Q(A)} < \epsilon} 
\geq - I_A(u) \, .
\]
Let $\eta > 0$ be such that $d - I_A(u) > \eta$. For $r$ sufficiently large
and $\tau = e^{r^d}$,
we have then 
\[
\begin{array}{ll}
\P \cro{\forall k \in K,  
\nor{ \pd * \xir - \pd * u_k}_{\infty, Q_k(A_r)}
\geq \epsilon  }
& \leq 
 \pare{1-e^{-r^d (I_A(u)+\eta)}}^{|\tilde{K}|} 
\\
& \approx e^{- \pare{\frac{1}{2 A_r}}^d 
\frac{\tau^{d-I_A(u)-\eta}}{(\log(\tau))^d}} \, .
\end{array} 
\]
Taking $\tau_n = \beta^n$ for some $\beta  >1$, the result follows from
Borel Cantelli lemma. 
\hspace*{\fill} \rule{2mm}{2mm}

\vspace{.5cm}
\noindent 
{\bf Step 2. A first lower bound.}

\begin{lemma}
\label{lbq.J1}
Let us define for $y \in \R$
\begin{equation} 
\label{J1-def}
\JJ_1(y) \triangleq 
\inf_{\mu \in \MM^{0}_1(\R^d)} 
\inf_{u \in \BB_1(\R^d)} \, 
\acc{ \LL(\mu) \, :
\bra{\mu,u}=y \, , \, I(u) <d }.
\end{equation} 
Let $\beta > 1$, and as before
$\tau_n = e^{r_n^d} = \beta^n$. Then, $\P$-a.s., $\forall y \in \R$,
$\forall \epsilon  > 0$,
\[
\varliminf_{\delta \rightarrow 0; \delta \in \Q} \, 
\varliminf_{n \rightarrow \infty}
\frac{1}{\tau_n} \log 
P_0 \cro{\va{\bra{L_{\tau_n},\pd * \bar{\xi}_{r_n}}-y} \leq \epsilon}
\geq - \JJ_1(y) \, .
\] 
Moreover, let $\JJ_1^{**}$ be the double Legendre transform of 
$\JJ_1$, then $\JJ_1^{**} = J$.
\end{lemma}
{\bf Proof.}
Let $\beta > 1$, $\delta > 0$, $A > 0$ 
and fix $u$ such that $I_A(u) < d$. Let 
$k$ be the index of the box of size $A_{r_n}$ associated by 
Lemma~\ref{a.s.comp.champ.lem} to $\delta$ and $\epsilon/4$. The center
of this box is denoted by $x_k=kA_{r_n}$. $\theta$ will
denote the shift on the Brownian trajectories, and $\sigma(D)$ will
denote the exit time of $D$.
\begin{eqnarray*}
 \lefteqn{P_0 \cro{ \va{\bra{L_{\tau},\pd*\xir} - y} \leq \epsilon}} 
\\
& & \geq P_0  \left[ \va{B_{\frac{\tau}{\log(\tau)}}-x_k}\leq 1;
\va{\frac{1}{\tau} \int_0^{\frac{\tau}{\log(\tau)}} \pd*\xir(B_s) \, ds} 
\leq \frac{\epsilon}{2}; \right.
 \\
& & 
\hspace*{2cm} \left .\sigma(Q_k(A_r)) \circ 
 \theta_{\frac{\tau}{\log(\tau)}} > \tau;
  \va{\frac{1}{\tau}\int_{\frac{\tau}{\log(\tau)}}^{\tau} \pd*\xir(B_s) \, ds - y} 
\leq \frac{\epsilon}{2} \right]  \, .
\end{eqnarray*}
Applying the Markov property at time $\tau/\log(\tau)$ yields,
\begin{eqnarray} 
\lefteqn{P_0 \cro{ \va{\bra{L_{\tau},\pd*\xir} - y} \leq \epsilon}}
\nonumber \\
& & \geq P_0 \cro{\va{B_{\frac{\tau}{\log(\tau)}}-x_k}\leq 1;
\va{\frac{1}{\tau} \int_0^{\frac{\tau}{\log(\tau)}} \pd*\xir(B_s) \, ds} 
\leq \frac{\epsilon}{2}}
\nonumber 
\\
& & \inf_{\va{x-x_k}\leq 1} \hspace{-.2cm}  P_x \cro{
\sigma(Q_k(A_r)) \geq \tau(1-\frac{1}{\log(\tau)}) ;
\va{\frac{1}{\tau}\int_{0}^{\tau(1-\frac{1}{\log(\tau)})} 
\hspace{-1cm} \pd*\xir(B_s) \, ds - y} 
\leq \frac{\epsilon}{2}} \, .
\label{lb1.produit}
\end{eqnarray}
Now, on $\sigma(Q_k(A_r)) \geq \tau(1-\frac{1}{\log(\tau)})$, 
\[
\va{\frac{1}{\tau}\int_{0}^{\tau(1-\frac{1}{\log(\tau)})}
\pd*\xir(B_s) \, ds
- \frac{1}{\tau}\int_{0}^{\tau(1-\frac{1}{\log(\tau)})}
\pd * u_k (B_s) \, ds}  \leq 
(1-\frac{1}{\log(\tau)}) \frac{\epsilon}{4} \, .
\]
Thus, for $\tau$ sufficiently large ($\frac{1}{\log(\tau)} 
< \frac{\epsilon}{8}$),
\begin{eqnarray*}
\lefteqn{\inf_{\va{x-x_k} \leq 1} P_x \cro{
\sigma(Q_k(A_r)) \geq \tau(1-\frac{1}{\log(\tau)}) ;
\va{\frac{1}{\tau}\int_{0}^{\tau(1-\frac{1}{\log(\tau)})} 
\hspace{-1cm} \pd*\xir(B_s) \, ds - y} 
\leq \frac{\epsilon}{2}}  }
\\
& & \geq 
\inf_{\va{x-x_k}\leq 1} P_x \cro{\sigma(Q_k(A_r)) 
\geq \tau(1-\frac{1}{\log(\tau)}) ;
\va{\frac{1}{\tau}\int_{0}^{\tau(1-\frac{1}{\log(\tau)})} 
\hspace{-1cm} \pd*u_k(B_s) \, ds - y} 
\leq \frac{\epsilon}{4}}
\\
&& 
\geq  \inf_{\va{x}\leq 1} P_x \cro{\sigma(Q(A_r)) 
\geq \tau(1-\frac{1}{\log(\tau)}) ;
\va{\bra{L_{\tau(1-\frac{1}{\log(\tau)})}; 
\pd*u } - y} 
\leq \frac{\epsilon}{8} } \, .
\end{eqnarray*} 
By the LDP lower bound 
for the Brownian occupation measure, we get then
\begin{equation} 
\begin{array}{l}
\displaystyle{\varliminf_{\tau \rightarrow \infty} \frac{1}{\tau}}
\log \inf_{\va{x-x_k}\leq 1} P_x \cro{
\sigma(Q_k(A_r)) \geq \tau(1-\frac{1}{\log(\tau)}) ;
\va{\frac{1}{\tau}\int_{0}^{\tau(1-\frac{1}{\log(\tau)})} 
\hspace{-1cm} \pd*\xir(B_s) \, ds - y} 
\leq \frac{\epsilon}{2}} 
\\[.3cm]
\hspace*{2cm}
 \geq - \displaystyle{\inf_{\mu \in \MM^0_1(Q(A))}}  \acc{ \LL(\mu):
\va{\bra{\mu, \pd *u}-y}< \frac{\epsilon}{8}} \, .
\end{array}
\label{lb1.produit.t2}
\end{equation} 
For the other term in \refeq{lb1.produit}, since
$\va{\frac{1}{\tau}\int_0^{\tau/\log(\tau)} \pd*\xir(B_s) \, ds}
\leq \frac{1}{\log(\tau)}$, we have for $\tau$ sufficiently large
($\frac{1}{\log(\tau)} \leq \frac{\epsilon}{2}$),
\begin{eqnarray}
\lefteqn{P_0 \cro{ \va{B_{\tau/\log(\tau)}-x_k}\leq 1;
\va{\frac{1}{\tau} \int_0^{\tau/\log(\tau)} \pd*\xir(B_s) \, ds}
\leq \frac{\epsilon}{2}}} 
\nonumber
\\
& & = 
P_0\cro{ \va{B_{\tau/\log(\tau)}-x_k}\leq 1}
\nonumber \\
& & = \int_{\nor{y-x_k}\leq 1} \exp\pare{- \frac{\nor{y}^2}{2 \tau/\log(\tau)}}
\, \frac{dy}{(2 \pi \tau/\log(\tau))^{d/2}} 
\nonumber \\
& & \geq 
 \frac{C(d)}{(2 \pi \tau/\log(\tau))^{d/2}}  
  \exp\pare{- \frac{(\tau/\log(\tau) +1)^2}{2 \tau/\log(\tau)}} \, .
\label{lb1.produit.t1}
\end{eqnarray} 
Putting \refeq{lb1.produit}, \refeq{lb1.produit.t2}, 
\refeq{lb1.produit.t1} together, we have that for any
$A > 0$, $u \in \BB_1(A)$ with $I_A(u) < d$,
$\delta >0$, $\epsilon >0$,  $\P$-a.s., for any $y \in \R$,  
\begin{equation} 
\label{lb1.u.fixe}
\varliminf_{\tau \rightarrow \infty} \frac{1}{\tau}
\log P_0 \cro{ \va{\bra{L_{\tau}, \pd*\xir}-y} \leq \epsilon}
\geq -\inf_{\mu \in \MM^0_1(Q(A))} \acc{ \LL(\mu): \va{\bra{\mu,\pd*u}-y} < 
\frac{\epsilon }{8}}
\, .
\end{equation} 
We would like to take the supremum over $u$ in the preceding expression.
Since the ``$\P$-a.s.'' depends on $u$, we have to restrict ourselves
to a countable subset of $\BB_1(A)$. 
\begin{lemma}
\label{sup.sur.u.lem}
For any $A > 0$, there exists a countable subset $\DD$ of 
$\BB_1(A)$, such that for any $u \in \BB_1(A)$, there is a 
sequence $(u_n)$ in $\DD$ satisfying
\begin{enumerate}
\item $\varlimsup_{n \rightarrow \infty} I_A(u_n) \leq I_A(u)$.
\item $\forall \mu \in \MM^0_1(Q(A))$, $\forall \delta  >  0$,
$\lim_{n \rightarrow  \infty} \bra{\mu, \pd * u_n} = \bra{\mu, \pd*u}$.
\end{enumerate}
\end{lemma} 
The proof of this lemma is given in section \ref{lemmes}. Lemma
\ref{sup.sur.u.lem} implies that $\forall \mu \in \MM^0_1(Q(A))$, and
$\forall \delta > 0$, $\forall \epsilon	> 0$, $\forall y \in \R$,
\[
\inf_{u \in \DD} \acc{ I_A(u): \va{\bra{\mu;\pd*u}-y} < \epsilon }
= \inf_{u \in \BB_1(A)} \acc{ I_A(u): \va{\bra{\mu;\pd*u}-y} < \epsilon }
\, .
\]
Thus, taking the supremum over $u \in \DD$ in \refeq{lb1.u.fixe},
we obtain that 
\begin{eqnarray*} 
\lefteqn{\varliminf_{\tau \rightarrow \infty} \frac{1}{\tau}
\log P_0 \cro{ \va{\bra{L_{\tau}, \pd*\xir}-y} \leq \epsilon}}
\\
& &
\geq -\inf_{\mu \in \MM^0_1(Q(A))} 
\acc{ \LL(\mu): \exists u \in \DD {\tt such that }
\va{\bra{\mu,\pd*u} -y} < 
\frac{\epsilon }{8} \mbox{ and } I_A(u) < d}
\\
& & \geq -\inf_{\mu \in \MM^0_1(Q(A))} \acc{ \LL(\mu): \,\, 
\inf_{u \in \BB_1(A)} \acc{ I_A(u): \va{\bra{\mu,\pd*u}-y} < 
\frac{\epsilon }{8}}
< d}
\end{eqnarray*} 

We now take $\delta$ to 0, $\delta \in \Q^+$. By Lemma~\ref{lem-tech3},
\begin{eqnarray*} 
\lefteqn{\varlimsup_{\delta \rightarrow 0, \delta \in \Q^+}\, 
\inf_{\mu \in \MM_1^0(Q(A))} 
\acc{  \LL(\mu)  \, :
\inf_{u \in \BB_1(A)} 
\acc{ I_A(u): \va{\bra{\mu,\psi_{\delta}*u}-y} < 
\frac{\epsilon}{8}} < d}}
\\
& & 
\leq 
\inf_{\mu \in \MM_1^0(Q(A))}
\acc{  \LL(\mu) \,  :
\inf_{u \in \BB_1(A)} 
\acc{ I_A(u): \va{\bra{\mu,u}-y} < \frac{\epsilon}{8}} < d}
\, .
\end{eqnarray*} 
We have thus proved that  $\forall A > 0$, $\forall \epsilon > 0$,
$\P$-a.s., 
$\forall y \in \R$, 
\begin{eqnarray*} 
\lefteqn{\varliminf_{\delta \rightarrow 0, \delta \in \Q^+} \,
\varliminf_{\tau \rightarrow \infty} \, 
\frac{1}{\tau} \log
P_0 \cro{ \va{\bra{L_{\tau},\pd*\xir}-y} \leq \epsilon}}
\\
& & 
\geq
- \inf_{\mu \in \MM_1^0(Q(A))} \, \inf_{u \in \BB_1(A)}
\acc{ \LL(\mu) \, : \, I_A(u) < d \, , \, 
\va{\bra{\mu,u}-y} < \frac{\epsilon}{8}}
\end{eqnarray*} 

Taking $A$ to infinity, it is easy to see that 
\begin{eqnarray*} 
\lefteqn{\lim_{A \rightarrow \infty, A\in \Q} 
\inf_{\mu \in \MM_1^0(Q(A))} \, \inf_{u \in \BB_1(A)}
\acc{ \LL(\mu) \, : \, I_A(u) < d \, , \, 
\va{\bra{\mu,u}-y} < \frac{\epsilon}{8}}}
\\
& & = 
\inf_{\mu \in \MM_1^0(\R^d)} \, \inf_{u \in \BB_1(\R^d)}
\acc{ \LL(\mu) \, : \, I(u) < d \, , \, 
\va{\bra{\mu,u}-y} < \frac{\epsilon}{8}} 
\leq \JJ_1(y) \, .
\end{eqnarray*}

We now prove that $\JJ_1^{**}(y)=J(y)$. Since $J=(-l)^*$,
it is enough to prove that $\JJ_1^*=-l$. It follows from the
large deviations estimates that $-\JJ_1 \leq -J$, so that 
$\JJ_1^* \leq J^*=(-l)^{**}=-l$, since $-l$ is convex continuous
(cf. Lemma~\ref{prop.l}). 
Hence, it remains to prove that $-l \leq \JJ_1^*$. 
A direct computation yields
\[
\JJ_1^*(\alpha) = - \inf_{u \in \BB_1}  \acc{\lambda(\alpha u,\R^d): 
\, I(u) < d} \, ,
\]
which is almost $-l$, except for the strict inequality ``$I(u) < d$'',
which we treat now.

Let $\alpha \in \R$
and $t$ be such that $t < -l(\alpha)$. It follows from the definition
of $l$ that $\exists u \in \BB_1$, $I(u) \leq d$,  and $\mu 
\in  \MM_1^0(\R^d)$, such that
$\LL(\mu)  - \alpha \bra{\mu,u} < -t$.
Let $A > 0$ 
and let us 
consider $d \mu_A(x)  \triangleq  A^{d} \frac{d \mu}{dx}(A x) \, dx$, and 
$u_A = u(A \cdot) $. Then $\mu_A \in \MM^0_1(\R^d)$,
$u_A \in \BB_1$, $\bra{\mu_A,u_A}=\bra{\mu,u}$,
$\LL(\mu_A) = A^2   \LL(\mu)$, and
$I(u_A) = A^{-d} I(u) \leq A^{-d} d < d$ for any $A >1$. Hence,
for any $A > 1$,
\begin{equation} 
- \JJ_1^*(\alpha) \leq 
\LL(\mu_A) - \alpha \bra{\mu_A,u_A} 
= A^2 \LL(\mu) - \alpha \bra{\mu,u} \, .
\end{equation} 
Therefore, $- \JJ_1^*(\alpha) \leq \LL(\mu)  - \alpha \bra{\mu,u} < -t$.
Since $t$ can be chosen arbitrarily in $]-\infty; -l(\alpha)[$, the result
follows.

\hspace*{\fill} \rule{2mm}{2mm} 

\noindent
{\bf Step 3. A sequence of lower bounds}.\\
For $y \in \R$, and $p \in \N^*$, let us define
\begin{eqnarray*} 
\lefteqn{ \DD_p(y) \triangleq 
\left\{ (\vec{\alpha},\vec{u},\vec{\mu}) \in [0,1]^p \times \BB_1^p 
\times \MM^0_1(\R^d)^p : \right. }
\\ 
& & \hspace*{3cm} \left. \sum_{j=1}^{p} \alpha_j = 1 \, , \, 
\sum_{j=1}^{p} \alpha_j \bra{\mu_j,u_j} =y \, , \,
\forall j, I(u_j) < d \right\} \, ,
\end{eqnarray*} 
\[
J_p(\vec{\alpha}, \vec{\mu}) \triangleq 
\sum_{j=1}^{p}  \alpha_j \LL(\mu_j) \, \, ,
\quad {\tt and} \quad
\JJ_p(y) \triangleq  \inf_{(\vec{\alpha},\vec{u}, \vec{\mu}) \in \DD_p(y)} 
J_p(\vec{\alpha}, \vec{\mu}) \, .
\]
\begin{lemma}
\label{lbqn.lem}
Let $\beta  > 1$, and let us define $\tau_n$ and $r_n$ by 
$\tau_n = e^{r_n^d} = \beta^n$. Then, $\forall \epsilon  > 0$,
$\P$-a.s., $\forall p \in \N$,
$\forall y \in \R$,
\[
\varliminf_{\delta \rightarrow 0; \delta \in \Q} \, 
\varliminf_{n \rightarrow \infty}
\frac{1}{\tau_n} \log 
P_0 \cro{\va{\bra{L_{\tau_n},\pd * \bar{\xi}_{r_n}}-y} \leq \epsilon}
\geq - \JJ_p(y) \, .
\] 
\end{lemma}
{\bf Proof}. The proof follows the same lines as step 2. 
Let $\beta  >1$, $\delta > 0$, $\epsilon > 0$, $A > 0$, 
$(\vec{\alpha},\vec{y})  \in 
]0,1[^p\times \R^p$, $\sum \alpha_j=1$, $\sum \alpha_j  y_j=y$, 
and $u_1, \cdots ,u_p$ such that $I_A(u_i) < d$ be fixed. Let 
$k_i$ be the indices of  the boxes  of size $A_{r_n}$ associated by 
Lemma~\ref{a.s.comp.champ.lem} to the $u_i$, $\delta$ and $\epsilon/6$. We 
divide the time interval $[0,\tau]$ in $p$ time intervals 
$[\tau_{i-1},\tau_i[$, where $\tau_i =\sum_{j=1}^{i} \alpha_j \tau$. In
the $i$-th time interval, we force the Brownian motion to
go fast (i.e., in time of order $\Delta_i = \alpha_i \tau/\log(\tau)$),
 from a neighborhood of $0$ to a neighborhood of $k_iA_r$, to remain
in $Q_{k_i}(A_r)$ during $\alpha_i \tau - 2\Delta_i$, and to return
in a neighborhood of 0 in time $\Delta_i$. We have then
\[
\begin{array}{l}
P_0 \cro{|\bra{L_{\tau};\pd * \xir}-y| < \epsilon} 
\\
\geq P_0 \left[ 
\forall i \in \{1,\cdots,p\}\, , \,\,
 |B_{\tau_{i-1}}| \leq 1; 
 \va{\frac{1}{\tau}\int_{\tau_{i-1}}
	^{\tau_{i-1} + \Delta_i} 
		\pd*\xir(B_s)	 \, ds }
			< \frac{\epsilon}{3p}; \right.
\\ \hspace*{2.5cm}
|B_{\tau_{i-1} + \Delta_i} - k_i A_r| \leq 1;
\sigma(Q_{k_i}(A_r))\circ \theta_{\tau_{i-1} + \Delta_i}
    > \tau_i  - \Delta_i	;
\\ \hspace*{2.5cm} 
 \va{\frac{1}{\tau} \int_{\tau_{i-1}+ \Delta_i}^
	     {\tau_i  - \Delta_i}  \pd * u_i(B_s)  ds - \alpha_i y_i}
	      < \frac{\epsilon}{6p};
\\ \hspace*{2.5cm} \left.
\va{\frac{1}{\tau}\int_{\tau_{i} - \Delta_i}
	^{\tau_{i} } \pd *\xir (B_s) \, ds }
			< \frac{\epsilon}{3p}
\right]
\end{array}
\]
By the Markov property and translation invariance, we get then
\[
P_0 \cro{|\bra{L_{\tau};\pd * \xir}-y| < \epsilon}
\geq \prod_{i=1}^{p} U_i V_i W_i \, ,
\]
where
\[\begin{array}[t]{l}
U_i = \inf_{|z| \leq 1}
	P_z \cro{\va{\frac{1}{\tau} \int_0^{\Delta_i} 
		 \pd *\xir (B_s) \, ds }
			< \frac{\epsilon}{3p};
	\va{B_{\Delta_i}-k_i A_r} \leq 1 }  \, ,
\\ 
V_i  = \inf_{z \in Q(1)}
	P_z \cro{ \sigma(Q(A_r)) 
		> \alpha_i \tau - 2 \Delta_i;
	\va{ \frac{1}{\tau} \int_0^{\alpha_i \tau - 2 \Delta_i} 
		 \pd* u_i(B_s) \, ds - \alpha_i y_i}
	         < \frac{\epsilon}{6p}} \, ,
\\
W_i  = \inf_{z \in Q_{k_i}(A_r)}
	P_z \cro{ 
	\va{\frac{1}{\tau} \int_0^{\Delta_i}
		\pd*\xir(B_s) \, ds }
			< \frac{\epsilon}{3p};
	\va{B_{\Delta_i}} \leq 1} \, .
\end{array}
\]
Exactly as in step 2, we can prove that 
\[
\varliminf_{\tau \rightarrow \infty} \frac{1}{\tau} \log(U_i) \geq 0 \, ,
\quad
\varliminf_{\tau \rightarrow \infty} \frac{1}{\tau} \log(W_i) \geq 0 \, ,
\]
\[
\varliminf_{\tau \rightarrow \infty} \frac{1}{\tau} \log(V_i)
\geq - \inf_{\mu_i \in \MM^0_1(Q(A))} \acc{\alpha_i \LL(\mu_i) :
\alpha_i \va{\bra{\mu_i, \pd*u_i}-y_i} < \frac{\epsilon}{6p}} \, .
\]
Taking now the supremum over $u_i \in \DD$, we get 
\begin{eqnarray*} 
\lefteqn{\varliminf_{\tau \rightarrow \infty} 
\frac{1}{\tau} \log   
P_0 \cro{|\bra{L_{\tau};\pd * \xir}-y| < \epsilon} }
\\
&& \geq - \inf_{\vec{\mu} \in \MM^0_1(Q(A))^p} \acc{
J_p(\vec{\alpha},\vec{\mu})\, : 
\, 
\inf_{u_i \in \DD_i} \, \max_i I_A(u_i) < d}
\, ,
\end{eqnarray*} 
where we have denoted $\DD_i= \acc{ u \in \DD; 
\alpha_i \va{\bra{\mu_i, \pd*u}-y_i} < \frac{\epsilon}{6p} }$.
Since
\[ \inf_{u_i \in \DD_i} 
\max_i I_A(u_i) = \max_i \inf_{u \in \DD_i} I_A(u) \, , 
\]
 the above infimum becomes by Lemma~\ref{sup.sur.u.lem} 
\[- \inf_{\vec{\mu} \in \MM^0_1(Q(A))^p} \,
\inf_{\vec{u} \in \BB_1(A)^p} 
\acc{J_p(\vec{\alpha},\vec{\mu}) \, :
\alpha_i \va{\bra{\mu_i, \pd*u_i}-y_i} < 
\frac{\epsilon}{6p}; I_A(u_i) < d} \, . 
\]
Taking $\delta$ to 0, we obtain
\begin{eqnarray*}
\lefteqn{- \inf_{\vec{\mu}, \vec{u}} \,
\acc{ J_p(\vec{\alpha},\vec{\mu})  \, : \, 
\alpha_i \va{\bra{\mu_i, u_i}-y_i} < \frac{\epsilon}{6p}\, ,\, 
I_A(u_i) < d} }
\\
&& \geq 
- \inf_{\vec{\mu} \in \MM^0_1(Q(A))^p} \,
\inf_{\vec{u} \in \BB_1(A)^p}  \acc{ J_p(\vec{\alpha},\vec{\mu}) \, : 
   \bra{\mu_i, u_i}=y_i \, ,   I_A(u_i) < d}  \, .
\end{eqnarray*} 
Optimizing  in $(\vec{\alpha},\vec{y})$, we are led to
\[
- \inf_{\vec{\alpha} \in ]0,1[^p} \, \inf_{\vec{\mu} \in \MM^0_1(Q(A))^p} \,
\inf_{\vec{u} \in \BB_1(A)^p}  \acc{J_p( \vec{\alpha},\vec{\mu})  \, : \,
\sum_{i=1}^p \alpha_i =1 \, , \,
\sum_{i=1}^{p} \alpha_i \bra{\mu_i, u_i} =y \, ,\, 
  I_A(u_i) < d} \,  .
\]
The proof of Lemma~\ref{lbqn.lem} follows after taking $A$ to infinity, and
noting that the infimum over $\alpha \in ]0,1[^d$ is the same as taking
$\alpha \in [0,1]^d$.
\hspace*{\fill} \rule{2mm}{2mm} 

\vspace{.5cm}
\noindent
{\bf Step 4. Conclusion}. \\
From Lemma~\ref{lbqn.lem}, Lemma~\ref{reg-lem}, and 
\refeq{rest.suite}, it is straightforward
to see that if $r(t)$ is defined as in Theorem~\ref{ldq.th}, 
then we have $\P$-a.s, for any $y \in \R$, $\epsilon > 0$, $p \in \N$,

\begin{equation}
\label{lbqp.gen}
\varliminf_{t \rightarrow \infty} \frac{r^2(t)}{t} 
\log P_0 \cro{ \va{\bra{L_t,\xi} - y   } <  \epsilon} 
\geq - \JJ_p(y) \, .
\end{equation} 
We now take $p$ to $\infty$.
\begin{lemma}. 
\label{prop.Jp}
\begin{enumerate}	
\item $\forall p \in \N$, $\forall y \in \R$,
$\JJ_{p+1}(y) \leq \JJ_p(y)$.
\item $\forall p \in \N$, $\forall \alpha \in [0,1]$, 
$\forall y_1,y_2 \in \R$,
\begin{equation}
\label{Jp.conv}
\JJ_{2p}(\alpha y_1 + (1-\alpha) y_2) 
\leq \alpha \JJ_p(y_1) + (1-\alpha) \JJ_p(y_2) \, .
\end{equation}	
\item Let $\JJ(y) \triangleq \lim_{p \rightarrow \infty}\  \JJ_p(y)$,
and $\tilde{\JJ}(y) \triangleq 
 \sup_{\epsilon > 0} \inf_{|z-y| \leq
\epsilon} \JJ(z)$ the greater lower 
semi-continuous minorant of $\JJ$. Then $\tilde{\JJ}=J$.
\item $\P$-a.s., $\forall y \in \R$, 
\[
\lim_{\epsilon \rightarrow 0} \,
\varliminf_{t \rightarrow \infty}
\frac{r^2(t)}{t} \log P_0 \cro{ \va{\bra{L_t,\xi}-y} \leq \epsilon }
\geq - J(y) \, .
\]
\end{enumerate} 
\end{lemma} 

\noindent
{\bf Proof of 1.}
For any $(\vec{\alpha}, \vec{u}, \vec{\mu}) \in \DD_p(y)$, 
and any $\nu$ with $\LL(\nu) < \infty$, we set
$\vec{\beta}\triangleq (\vec{\alpha}, 0)$, 
$\vec{w}\triangleq (\vec{u},0)$ and $\vec{\pi} \triangleq (\vec{\mu},\nu)$.
We note that $(\vec{\beta},\vec{w},\vec{\pi}) \in \DD_{p+1}(y)$. Thus,
\[ 
\JJ_{p+1}(y) \leq J_{p+1}(\vec{\beta},\vec{\pi}) = 
J_p(\vec{\alpha}, \vec{\mu}).
\]
Taking the infimum over $\DD_p(y)$ yields $\JJ_{p+1}(y) \leq  \JJ_p(y)$.

\vspace{.5cm}
\noindent
{\bf Proof of 2.}
In the same way, let $\alpha \in [0,1]$ and $y_1,y_2 \in \R$ be
fixed. For any $(\vec{\beta_1},\vec{u_1},\vec{\mu_1}) \in \DD_p(y_1)$, and
any $(\vec{\beta_2 },\vec{u_2},\vec{\mu_2}) \in \DD_p(y_2)$, we set
$\vec{\lambda} \triangleq (\alpha \vec{\beta_1}, (1-\alpha) \vec{\beta_2})$,
$\vec{v} \triangleq  (\vec{u_1},\vec{u_2})$, and
$\vec{\nu} \triangleq  (\vec{\mu_1},\vec{\mu_2})$. Note that 
 $(\vec{\lambda},\vec{v},\vec{\nu})
\in \DD_{2p}(\alpha y_1+(1-\alpha)y_2)$. Thus,
\[
\JJ_{2p}(\alpha y_1+(1-\alpha)y_2) 
\leq J_{2p}(\vec{\lambda},\vec{\nu}) 
= \alpha J_p(\vec{\beta_1},\vec{\mu_1}) + 
(1-\alpha) J_p(\vec{\beta_2},\vec{\mu_2})
\,.
\]
Taking the infimum over elements of $\DD_p(y_1)$ and $\DD_p(y_2)$, leads
to \refeq{Jp.conv}.

\vspace{.5cm}
\noindent
{\bf Proof of 3 and 4.}
Taking $p$ to $\infty$ in \refeq{lbqp.gen} yields that $\P$-a.s,  for any
$y \in \R$ and $\epsilon > 0$, 
\begin{equation}
\varliminf_{t \rightarrow \infty} \frac{r^2(t)}{t} 
\log P_0 \cro{ \va{\bra{L_t,\xi} - y   } \leq \epsilon} 
\geq - \JJ(y) \, .
\end{equation} 
From this, it follows easily that $\P$-a.s, $\forall y \in \R$, 
\begin{equation}
\label{lbqJ}
\lim_{\epsilon \rightarrow 0} \, 
\varliminf_{t \rightarrow \infty} \frac{r^2(t)}{t} 
\log P_0 \cro{ \va{\bra{L_t,\xi} - y   } \leq \epsilon}
\geq - \tilde{\JJ}(y) 
\end{equation} 
 From the large deviations upper bound, we have then 
$\tilde{\JJ} \geq J$. On the other hand, we also have 
$\tilde{\JJ} \leq \JJ \leq \JJ_1$, so that 
$\tilde{\JJ}^{**} \leq \JJ_1^{**} = J(y)$ by Lemma~\ref{lbq.J1}. 
Taking $p$ to $\infty$ in \refeq{Jp.conv}, we obtain that
$\JJ$ is convex. Thus, $\tilde{\JJ}$ is convex and 
lower semi-continuous. Thus, $\tilde{\JJ}=\tilde{\JJ}^{**}$, and 
$\tilde{\JJ} \leq J$. Finally, $\tilde{\JJ} =  J$.
\hspace*{\fill} \rule{2mm}{2mm}


\section{Technical Lemmas}
\label{lemmes}

\noindent
{\bf Proof of Lemma~\ref{lem-tech3}.} \\
We can assume that $\LL(\mu)<\infty$. Let then $\varphi=\sqrt{d\mu/dx}$, 
and as $u$ is bounded by 1,
\begin{equation} 
\label{eq1.6}
|\int\pare{\pd*u-u}d\mu|\le \int_{\R^d} |\pd*\varphi^2(x)-\varphi^2(x)|dx.
\end{equation} 
Now, for any $\epsilon_1>0$, and any $x\in \R^d$
\begin{align}\label{eq1.7}
\int \pd(y)|\varphi^2(x-y)-\varphi^2(x)|dy\le \frac{1}{2\epsilon_1}&
\int \pd(y)\pare{\varphi(x-y)-\varphi(x)}^2dy\notag\\
&+\epsilon_1\pare{\pd*\varphi^2(x)+\varphi^2(x)}.
\end{align}
Also,
\begin{equation} 
\label{eq1.8}
|\varphi(x-y)-\varphi(x)|=|\int_0^1 \nabla \varphi(x-ty).y dt|\le ||y||
\pare{\int_0^1||\nabla \varphi(x-ty)||^2dt}^{1/2}.
\end{equation} 
Thus,
\begin{align}
\label{eq1.9}
\int_{\R^d} |\pd*\varphi^2(x)-\varphi^2(x)|dx&\le
\frac{1}{2\epsilon_1}
\int_0^1dt\int_{\R^d}dy \pd(y)||y||^2\int_{\R^d}dx ||\nabla \varphi(x-ty)||^2
+2\epsilon_1\notag\\
&\le \frac{\int \pd(y)||y||^2dy}{2\epsilon_1}\int_{\R^d}||\nabla \varphi(x)||^2dx+
2\epsilon_1.
\end{align}
There is a constant $c_0$ such that $\int\pd(y)||y||^2dy=c_0 \delta ^2$,
and the result follows.
 \hspace*{\fill} \rule{2mm}{2mm}

\vspace{.5cm}
\noindent
{\bf Proof of Lemma~\ref{cont-bra}.} \\
Let $(u_n)$ be a sequence converging weakly to $u \in \BB_1(A)$, and
$(\mu_n)$ a sequence converging weakly to $\mu \in \MM_1(Q(A))$. 
We think of $u_n$ and $u$ as vanishing outside $Q(A)$.
For any $\delta >0$, $(\pd * u_n)_n$ is an equicontinuous, uniformly
bounded sequence converging pointwise to $\pd*u$. By
Ascoli-Arzel\`a, we have
\[
\lim_{n\to\infty} \nor{\pd*u_n -\pd*u}_{\infty, Q(A)}=0.
\]
The result follows then from the inequality
\[
|\bra{\mu_n, \pd * u_n} - \bra{\mu, \pd*u}|
\leq \nor{ \pd * u_n- \pd*u}_{\infty, Q(A)} 
+ |\bra{\mu_n, \pd*u} - \bra{\mu,\pd*u}| \, .
\]
\hspace*{\fill} \rule{2mm}{2mm}

\vspace{.5cm}
\noindent
{\bf Proof of Lemma~\ref{pb-dual}.} \\
Set 
\[
\II_1(y) = \sup_{a > 0} \acc{ ay - \int_{\TT(A)} \Lambda(a f(x)) \, dx } \, ,
\]
\[
\II_2(y) = \inf_{u \in \BB_1(A)} \acc{ I_A(u): \int_{Q(A)} f(x) u(x) \, dx 
\geq y} \, .
\]
Note that
\[
\II_2(y)  = \inf_{u \in \BB_1(A)} \, \sup_{a > 0}
\acc{ I_A(u) + a \pare{y - \int_{Q(A)}  f(x) u(x) \, dx }}    \, .
\]
Inverting the infimum and the supremum in the preceding expression, we
obtain  $\II_1(y)$. Hence $\II_2(y) \geq \II_1(y)$, and $\II_2(y)$ and
$\II_1(y)$ are dual optimization problems. 

Since $\forall a > 0$, $\Lambda(a) \leq a M$, it follows from
the definition of $\II_1$ that $\forall a > 0$,
$y \leq \frac{\II_1(y)}{a} + M$. Hence if $\II_1(y) < \infty$, then
$y \leq M$. In other words, for $y > M$, $+ \infty = \II_1(y) 
\leq \II_2(y)$.

 For $y < M$, note that $\II_2(y) \leq I_A(u \equiv y) = |Q(A)| H(y) < \infty$.
Moreover, the infimum in $\II_2$ is actually a minimum. Actually,
$\CC_y = \acc{ u \in \BB_1(A); \bra{f,u} \geq y}$ is compact 
in weak topology. Indeed, let 
$(u_n)$ be a sequence in $\CC_y$. It follows from Banach-Alaoglu theorem
that $u_n$ converges weakly to $u \in \BB_1(A)$. Hence $\bra{f,u_n}
\rightarrow \bra{f,u}$, and $u\in \CC_y$. $I_A$ being lower semi-continuous,
the infimum of $I_A$ on $\CC_y$ is a minimum, as soon as 
$\CC_y$ is not empty, which is actually the case, since $u \equiv y$
belongs to $\CC_y$. Moreover, for $y < M$, the Slater condition
(see for instance  Theorem~6.7 in \cite{jahn}) is satisfied
by $u\equiv z$, for $z \in ]y;M[$.   The identity between 
$\II_1(y)$ and $\II_2(y)$ follows then from standard results in convex
optimization. 

We have thus proved that $\II_1=\II_2$, except on $y=M$. But, note that 
$\II_1$ and $\II_2$ are obviously increasing on $\R^+$. $\II_1$ is
clearly lower semi-continuous, 
and the same is true for $\II_2$. Indeed, let $(y_n)$ 
a sequence converging to $y$, and let $L$ be such that 
$\liminf_{n \rightarrow  \infty} \II_2(y_n) < L$. 
We can then find a (sub)sequence
$(u_n)$ in $\BB_1(A)$, with $\bra{f, u_n} \geq y_n$, and $
I_A(u_n) < L$ for sufficiently large $n$.  $\BB_1(A)$ being 
weakly compact, there exists $u \in \BB_1(A)$ such that $u_n$ converges
weakly to $u$. Hence $\bra{f,u_n} \rightarrow  \bra{f,u}$, 
so that $u \in \CC_y$.
Therefore, $\II_2(y) \leq I_A(u) \leq 
\liminf_{n \rightarrow  \infty} I_A(u_n)$
by lower continuity of $I_A$. Hence $\II_2(y) \leq  L$ for any 
$L > \liminf_{n \rightarrow \infty} \II_2(y_n)$.

 By lower semi-continuity
and monotonicity, we have
that 
\[ \II_1(M) \leq \liminf_{y_n \nearrow M} \II_1(y_n)
  \leq \limsup_{y_n \nearrow M} \II_1(y_n)  \leq \II_1(M) \, ,
\]
and the same holds true for $\II_2$. Hence $\II_1(M)=\II_2(M)$.
\hspace*{\fill} \rule{2mm}{2mm}

\vspace{.5cm}
\noindent 
{\bf Proof of Lemma~\ref{lem2.2}.}  \\
We can assume that there is $L<\infty$ such that
\[
\limsup_{A\to\infty} \, \inf_{z \geq y} \bar{\II}_{\delta,A}(z) = L. 
\]
For sufficiently large $A$, let $
\mu_A,u_A\in \MM_1(\TT(A))\times\BB_1(Q(A))$ be such that
\begin{equation} 
\label{eq5.6} 
 I_{A}(u_A)+ \LL_{A}(\mu_A) < L + \frac{1}{A}
\quad {\tt and} \quad
\int_{Q(A)} u_A \, d\pd^A*\mu_A \geq y.
\end{equation} 
Note that changing $u_A$ on $\partial Q(A)$ does not change anything in the
above expres\-sion, and we can as well assume that $u_A \equiv 0$
on $\partial Q(A)$. We extend   $u_A$ outside $Q(A)$ by periodization.
 Following Lemma  3.5 of~\cite{donsker-varadhan}, it is possible to
translate both  $\mu_A$ and $u_A$ by the same amount
--we still call  $\mu_A,u_A$ the translates-- in such a way that
\begin{equation} 
\label{eq5.7} 
\mu_A(\partial_A Q(A))\le \frac{2d}{\sqrt{A}} \, ,
\end{equation} 
where $\partial_A Q(A)=\bigcup_{i=1}^d\{ -\frac{A}{2}\le x_i\le -\frac{A}{2}+
\sqrt{A}\}\cup \{ \frac{A}{2}-\sqrt{A}\le x_i\le \frac{A}{2}\}.$
And there is a measure $\tilde \mu_A$ with Dirichlet boundary on
$Q_0(A)$ such that
\begin{equation} 
\label{eq5.8}
|\LL_A(\mu_A)-\LL_A(\tilde \mu_A)|\le \frac{2d}{\sqrt{A}},\quad
{\tt and} \quad
\tilde \mu_A\big|_{Q(A)\backslash \partial_A Q(A)}
=\mu_A\big|_{Q(A)\backslash \partial_A Q(A)}.
\end{equation} 
We denote by $\tilde u_A$ the function vanishing on $\partial_A Q(A)$,
and equal to $u_A$ on $Q(A)\backslash \partial_A Q(A)$. 
Note that $I_A(\tilde u_A)\le I_A(u_A)$ and
\begin{equation} 
\label{eq5.9}
|\int \pd^A*u_A d\mu_A-\int \pd*\tilde u_A d\tilde \mu_A|\le
|\pd*u_A|_{\infty} \mu_A(\partial_A Q(A))\le \frac{2d}{\sqrt{A}}.
\end{equation} 
Thus, $\int \pd*\tilde u_A d\tilde \mu_A \geq y- \epsilon$ 
for $\epsilon>2d/\sqrt{A}$,
and $\tilde{\mu}_A,\tilde{u}_A\in \MM_1^0(\R^d) \times \BB_1$. 
This completes the proof.
\hspace*{\fill} \rule{2mm}{2mm} 

\vspace{.5cm}
\noindent 
{\bf Proof of Lemma~\ref{fr.minvp.lem}} 
\[
\P \cro{
\min_{1 \leq k \leq  \NN} 
\acc{\lambda \pare{\alpha \pd  *\xir,
 Q_k(A)}} \leq x}
\leq  \sum_{1 \leq k \leq  \NN}
\P \cro{\lambda \pare{\alpha \pd  *\xir,
 Q_k(A)} \leq x} \, .
\] 
  Note that if $A$ is a  multiple integer
of $1/r$, all the random variables appearing in the sum  have the same law.
Since $Q_k(A) \subset \frac{[rkA]}{r} + Q(\frac{\cro{rA}}{r}
+ \frac{2}{r})$, we have for any $r>2$, 
\[ 
\begin{array}{l} 
\P \cro{
\displaystyle{\min_{1 \leq k \leq \NN} }
\acc{\lambda \pare{\alpha \pd  *\xir,
 Q_k(A)}} \leq x}
\\[.2cm]
 \hspace*{2cm} \leq \pare{\frac{R\tau +A}{A}}^d  
\P \cro{\lambda \pare{\alpha \pd  *\xir,
 Q(\frac{\cro{Ar}}{r}
+ \frac{2}{r})} \leq x} 
\\[.2cm]
 \hspace*{2cm}
 \leq \pare{\frac{R\tau +A}{A}}^d \P \cro{\lambda \pare{\alpha \pd  *\xir,Q(A+1)} \leq x}
\end{array} 
\]
If we impose the relation $\tau=\exp(r^d)$, the result follows from 
the LDP for $\xir$ (lemma \ref{lem.champ-LDP}), since by 
Lemma~\ref{cor-tech1}, $\acc{u \in \BB_1(A+1), 
\lambda \pare{\alpha \pd  * u,Q(A+1)} \leq x}$ is closed. 
\hspace*{\fill} \rule{2mm}{2mm}

\vspace{.5cm}
\noindent
{\bf Proof of Lemma~\ref{sup.sur.u.lem}.} \\
For any $n \in \N$, let us partition $Q(A)$ into $N_n$ 
dyadic cubes of order $n$, denoted by $I^{(n)}_j$. Let 
\[
\DD_n \triangleq  \acc{ \sum_{j=1}^{N_n} \frac{k_j}{2^n} \ind_{I^{(n)}_j};
k_j \in [-2^n, 2^n] \cap \Z} \, , 
\mbox{ and }  
\DD \triangleq \cup_n \DD_n \, .
\]
For any real $x$, let $\lfloor x \rfloor$ be the nearest 
integer of $x$ in the interval $[0,x]$, i.e.
\[
\lfloor x \rfloor = \left\{ 
\begin{array}{ll} k \mbox{ if } 0 \leq k \leq x < k+1 \, , \\
		  k \mbox{ if } k-1 < x \leq k \leq 0 \, .
\end{array} 
\right.
\]
We associate to any $u \in \BB_1(A)$, the function $u_n$ of $\DD_n$
defined by 
\[ u_n \triangleq \sum_{j=1}^{N_n} \frac{\lfloor 2^n \bar{u}^{(n)}_j 
\rfloor}{2^n} \ind_{I^{(n)}_j} \, , \,
\mbox{ where }
\bar{u}^{(n)}_j \triangleq \frac{1}{|I^{(n)}_j|} 
\int_{I^{(n)}_j} u(x) \, dx \, .
\]
Since $H$ is convex, increasing on $\R^+$, decreasing on 
$\R^-$, we get 
\[ 
I_A(u_n)  = \sum_{j=1}^{N_n} H\pare{\frac{\lfloor 2^n \bar{u}^{(n)}_j 
\rfloor}{2^n}} |I^{(n)}_j|  
 \leq   \sum_{j=1}^{N_n} H( \bar{u}^{(n)}_j)  |I^{(n)}_j| 
 \leq  I_A(u) \, .
\]
Moreover, 
\[
\begin{array}{l}
\va{\bra{\pd*\mu,u_n}-\bra{\pd*\mu,u}}
\\
 \leq \va{ \displaystyle{\sum_{j=1}^{N_n}} 
(\frac{\lfloor 2^n \bar{u}^{(n)}_j 
\rfloor}{2^n} -  \bar{u}^{(n)}_j) \int_{I^{(n)}_j} \pd*\mu(x) \, dx }
+ \va{ \displaystyle{\sum_{j=1}^{N_n} }
\int_{I^{(n)}_j} \pd*\mu(x) (\bar{u}^{(n)}_j-u(x)) \, dx }
\\ 
 \leq \frac{1}{2^n} + \va{\displaystyle{\sum_{j=1}^{N_n}}
\int_{I^{(n)}_j} u(x) \int_{I^{(n)}_j} 
(\pd *\mu(y) - \pd*\mu(x)) \, \frac{dy}{|I^{(n)}_j|} \, dx }
\\
 \leq \frac{1}{2^n} + \omega(\pd*\mu, \frac{1}{2^n}) \int_{Q(A)} u(x) \, dx
\, ,
\end{array}
\]
where $\omega(\pd*\mu, \frac{1}{2^n})$ is the modulus of continuity
of $ \pd*\mu$ on $Q(A)$.
\hspace*{\fill} \rule{2mm}{2mm}

\vspace{.5cm}
\begin{lemma} 
\label{prop.l}
$l$ is concave and continuous, and takes negative values. Moreover, $l(0)=0$.
$J$ is convex and lower semi-continuous, increasing on $\R^+$, decreasing
on $\R^-$. $J(0)=0$, $J = \infty$ outside $[m,M]$, and $J$ is finite
on $]m,M[$.
\end{lemma}
{\bf Proof}. $l$ is concave as the infimum of affine functions. 
$I(0)=0 \leq d$. Hence,
$l(\alpha) \leq \lambda(0,\R^d)=0$. Moreover, $l(0)=\lambda(0,\R^d)=0$.
Since $H(y)=+ \infty$ for $y \notin [m,M]$, 
when $u$ is such that $I(u) \leq d$,
one also has $m \leq u(x) \leq M$
$dx$-a.s.. It follows easily that
\begin{equation} 
\label{borne.l.eq}
- \alpha M \leq l(\alpha) \leq 0 \mbox{ for } \alpha \geq 0 \, , 
\quad {\tt and} \quad
- \alpha m \leq l(\alpha) \leq 0 \mbox{ for } \alpha \leq 0 \, .
\end{equation} 
The continuity of $l$ is then a consequence of its concavity, and of
the fact that it is everywhere finite.

Now, $J$ is lower semi-continuous and convex as supremum of affine functions.
$J(y) \geq 0.y + l(0)=0$. $J(0)=\sup_{\alpha \in \R} l(\alpha)=0$.
One deduces from \refeq{borne.l.eq} that 
\[ 
\begin{array}{lcl}
J(y) & \geq & 
\max \acc{ \sup_{\alpha \leq 0}\acc{\alpha(y- m)} \, , 
\sup_{\alpha \geq 0} \acc{\alpha(y- M)}} 
\\
& \geq  & +\infty \mbox{ for } y \notin [m;M] 
\end{array}
\]
Let us prove that $J$ is finite on $]m,M[$. Since $J \leq \JJ_1$, it is 
enough to prove that $\JJ_1$ is finite on $]m,M[$. But for any 
$y \in ]m,M[$, $H(y) < \infty$, and one can find $\epsilon > 0$ such 
that $H(y) \epsilon^d < d$. Thus, $u \triangleq y \ind_{Q(\epsilon)}$ is
such that $I(u) < d$, so that $\JJ_1(y) \leq \LL(\mu)$ for any $\mu \in 
\MM^0_1(Q(\epsilon))$. 

It remains to prove the monotonicity of $J$. 
For $y \in \R^+$, $\sup_{\alpha \leq 0} \acc{\alpha  y +l(\alpha)}
\leq 0 \leq J(y)$. Hence $J(y)=\sup_{\alpha \geq 0}\acc{\alpha y
+ l(\alpha)}$, and $J$ is increasing on $\R^{+}$. 
\hspace*{\fill} \rule{2mm}{2mm}   

\vspace{.5cm}
\begin{lemma}.\\ 
\label{cor-tech1} $\forall A > 0$, $\forall \delta > 0$, $\forall x \in \R$,
$\forall \alpha \in \R$,
$\acc{u \in \BB_1(A); \lambda(\alpha \pd*u, Q(A)) \leq x}$ is 
compact in weak topology.
\end{lemma}
{\bf Proof}. Since $\BB_1(A)$ is weakly compact, it is enough to 
prove that $u \in \BB_1(A) \mapsto \lambda(\alpha \pd*u, Q(A))$ is
lower semi-continuous. 
Let then $(u_n)$ a sequence in $\BB_1(A)$ weakly converging
to $u$, and let $L > \liminf_{n \rightarrow \infty}
 \lambda(\alpha \pd*u_n, Q(A))$. By definition of 
$\lambda(\alpha \pd*u_n, Q(A))$, one can then find a (sub)sequence
of probability measures $(\mu_n) \in \MM^0_1(Q(A))$ such that for sufficiently
large $n$, $\LL(\mu_n) + \alpha \bra{\mu_n, \pd * u_n} < L$. For such
$n$, $\LL(\mu_n) \leq L + |\alpha|$, and there exists $\mu \in \MM^0_1(Q(A))$,
and a subsequence $(n_k)$ such that $\mu_{n_k}$ converges weakly to $\mu$.
It follows then from Lemma~\ref{cont-bra}, and the lower
semi-continuity of $\LL$  that
\[
\lambda(\alpha   \pd*u, Q(A)) \leq 
\LL(\mu) + \alpha \bra{\mu, \pd*u} 
\leq \liminf_{k \rightarrow \infty} 
(\LL(\mu_{n_k}) + \alpha \bra{\mu_{n_k}, \pd*u_{n_k}})
\leq  L \, .
\]
The proof is completed as $L$ tends to $\liminf_{n \rightarrow \infty} 
\lambda(\alpha   \pd*u_n, Q(A))$.
\hspace*{\fill} \rule{2mm}{2mm}

\end{document}